\documentclass[12pt]{amsart}

\usepackage{amssymb}
\usepackage{graphicx}
\usepackage{ifthen}
\usepackage{pict2e}
\usepackage{xargs}
\usepackage{xspace}

\newcommand{\definedterm}[1]{\emph{#1}}

\newcommand{\Bairespace}[1][]{
  \ifthenelse{\equal{#1}{}}{\functions{\N}{\N}}{\functions{#1}{\N}}
}
\newcommand{\Bairetree}[1][]{
  \ifthenelse{\equal{#1}{}}{\functions{<\N}{\N}}{\functions{#1}{\N}}
}
\newcommand{\ball}[3][]{\calB_{#1}(#2, #3)}

\newcommand{\calB}{\mathcal{B}}
\newcommand{\calN}{\mathcal{N}}

\newcommand{\calU}{\mathcal{U}}
\newcommand{\calV}{\mathcal{V}}
\newcommand{\Cantorspace}[1][]{
  \ifthenelse{\equal{#1}{}}{\functions{\N}{2}}{\functions{#1}{2}}
}
\newcommand{\Cantortree}[1][]{
  \ifthenelse{\equal{#1}{}}{\functions{<\N}{2}}{\functions{#1}{2}}
}

\newcommand{\closure}[1]{\overline{#1}}

\newcommand{\compactifiedBairespace}{\Bairespace_*}

\newcommand{\compactifiedextendedBairespace}{\extendedBairespace_*}
\newcommand{\composition}{\circ}

\newcommandx{\concatenation}[2][1 = undefined, 2 = undefined]{
  \ifthenelse{\equal{#1}{undefined}}{{}\smallfrown}{
    \ifthenelse{\equal{#2}{undefined}}{\bigoplus #1}{\bigoplus_{#1} #2}
  }
}
\newcommandx{\constant}[2][2 =]{\ifthenelse{\equal{#2}{}}{c_{#1}}{c_{#1, #2}}}
\newcommandx{\constantfunction}[3][2 =, 3 =]{
  \ifthenelse{\equal{#2}{}}{c \from #1 \to \image{c}{#1}}{c_{#3} \from #1 \to #2}
}

\newcommandx{\convolution}[2][1 = undefined, 2 = undefined]{
  \ifthenelse{\equal{#1}{undefined}}{\mathrel{*}}{
    \ifthenelse{\equal{#2}{undefined}}{\bigotimes #1}{\bigotimes_{#1} #2}
  }
}
\newcommandx{\Deltaclass}[2][1=,2=]{
  \ifthenelse{\equal{#2}{}}{\mathbf{\Delta}_{#1}}{\mathbf{\Delta}^{#1}_{#2}}
}
\newcommand{\diameter}[2][]{\mathrm{diam}_{#1} \thinspace #2}

\newcommand{\discrete}{p}

\newcommandx{\disjointunion}[2][1 =, 2 =]{
  \ifthenelse{\equal{#1}{}}{\sqcup}{
    \ifthenelse{\equal{#2}{}}{\bigsqcup #1}{{\bigsqcup_{#1} #2}}
  }
}

\newcommand{\domain}[1]{\mathrm{dom}(#1)}

\newcommand{\emptysequence}{\emptyset}

\newcommand{\extendedBairespace}{\Bairespace[\le \N]}

\newcommand{\extendedby}{\sqsubseteq}

\newcommand{\extendedextensions}[1]{\extensions{#1}^*}

\newcommand{\extension}[1]{\overline{#1}}
\newcommand{\extensions}[2][]{
  \ifthenelse{\equal{#1}{}}{\calN_{#2}}{\extendedextensions{#2} \intersection #1}
}

\newcommand{\from}{\colon}
\newcommand{\Fsigma}{$F_\sigma$\xspace}
\newcommandx{\functions}[3][3 =]{
  \ifthenelse{\equal{#3}{}}{#2^{#1}}{#2^{#1}_{#3}}
}

\newcommand{\Gdelta}{$G_\delta$\xspace}
\newcommand{\goesto}{\rightarrow}

\newcommand{\GzeroN}[1][]{
  \ifthenelse{\equal{#1}{}}{\mathbb{G}_0^\N}{\mathbb{G}_{0, #1}^{\N}}
}

\newcommand{\identity}[1]{\mathrm{id}_{#1}}
\newcommand{\image}[2]{#1(#2)}
\newcommand{\inclusion}[2]{\iota_{#1, #2}}

\newcommandx{\identityfunction}[2][2 =]{
  \ifthenelse{\equal{#2}{}}{\mathrm{id} \from #1 \to #1}{\mathrm{id} \from #1 \to #2}
}
\newcommand{\infimum}[2][]{
  \ifthenelse{\equal{#1}{}}{\inf #1}{\inf_{#1}{#2}}
}
\newcommand{\injections}[1]{\functions{\N}{(\N)}}
\newcommandx{\intersection}[2][1 =, 2 =]{
  \ifthenelse{\equal{#1}{}}{\cap}{
    \ifthenelse{\equal{#2}{}}{\bigcap #1}{{\bigcap_{#1} #2}}
  }
}
\newcommand{\inverse}[1]{#1^{-1}}

\newcommand{\length}[1]{|#1|}
\newcommandx{\limit}[2][1 =, 2 =]{
  \ifthenelse{\equal{#1}{}}{\lim}{
    \ifthenelse{\equal{#2}{}}{\lim #1}{{\lim_{#1} #2}}
  }
}

\newcommand{\mathand}{\text{ and }}

\newcommand{\meet}{\wedge}

\newcommand{\N}{\mathbb{N}}
\newcommand{\metric}[1]{d_{#1}}

\newcommand{\onepointcompactification}[1]{#1_*}

\newcommand{\pair}[2]{(#1, #2)}
\newcommandx{\Piclass}[2][1=,2=]{
  \ifthenelse{\equal{#2}{}}{\mathbf{\Pi}_{#1}}{\mathbf{\Pi}^{#1}_{#2}}
}

\newcommand{\preimage}[2]{#1^{-1}(#2)}
\newcommandx{\product}[2][1 =, 2 =]{
  \ifthenelse{\equal{#1}{}}{\times}{
    \ifthenelse{\equal{#2}{}}{\prod #1}{{\prod_{#1} #2}}
  }
}
\newcommandx{\projection}[2][1 =, 2 =]{
  \ifthenelse{\equal{#1}{}}{\mathrm{proj}}{
    \ifthenelse{\equal{#2}{}}{\projection_{#1}}{
      \image{\projection[#1]}{#2}
    }
  }
}
\renewcommandx{\restriction}[3][3 = undefined]{
  \ifthenelse{\equal{#3}{undefined}}{#1 \upharpoonright #2}{#1 \upharpoonright #2 \to #3}
}
\newcommandx{\sequence}[2][2 = undefined]{
  \ifthenelse{\equal{#2}{undefined}}{(#1)}{
    (#1)_{#2}
  }
}
\newcommandx{\set}[2][2 = undefined]{
  \ifthenelse{\equal{#2}{undefined}}{\{ #1 \}}{
    \{ #1 \suchthat #2 \}
  }
}
\newcommand{\setcomplement}[1]{\twiddle #1}
\newcommandx{\sets}[3][3 =]{
  \ifthenelse{\equal{#3}{}}{[#2]^{#1}}{[#2]^{#1}_{#3}}
}

\newcommandx{\Sigmaclass}[2][1=,2=]{
  \ifthenelse{\equal{#2}{}}{\mathbf{\Sigma}_{#1}}{\mathbf{\Sigma}^{#1}_{#2}}
}
\newcommand{\strictlyextendedby}{\sqsubset}

\newcommand{\suchthat}{\mid}

\newcommand{\twiddle}{\raisebox{1pt}{\scalebox{.75}{$\mathord{\sim}$}}}
\newcommandx{\union}[2][1 =, 2 =]{
  \ifthenelse{\equal{#1}{}}{\cup}{
    \ifthenelse{\equal{#2}{}}{\bigcup #1}{{\bigcup_{#1} #2}}
  }
}

\newcommand{\Baire}{Baire\xspace}
\newcommand{\Borel}{Bor\-el\xspace}
\newcommand{\Brouwer}{Brouw\-er\xspace}

\newcommand{\Hurewicz}{Hur\-e\-wicz\xspace}

\newcommand{\Jayne}{Jayne\xspace}

\newcommand{\Polish}{Po\-lish\xspace}
\newcommand{\Rogers}{Rog\-ers\xspace}

\newenvironment{lemmaproof}{
  
  \begin{proof}
}{\end{proof}}

\newenvironment{propositionproof}{
  
  \begin{proof}
}{\end{proof}}

\newenvironment{theoremproof}{
  
  \begin{proof}
}{\end{proof}}

\newtheorem{lemma}{Lemma}[section]
\newtheorem{proposition}[lemma]{Proposition}

\newtheorem{theorem}[lemma]{Theorem}

\newtheorem{introtheorem}{Theorem}

\theoremstyle{definition}
\newtheorem{remark}[lemma]{Remark}

\begin{document}

%% Front matter

\begin{abstract}
  We provide a finite basis for the class of \Borel functions that are
  not in the first \Baire class, as well as the class of \Borel functions
  that are not $\sigma$-continuous with closed witnesses.
\end{abstract}

\author[R. Carroy]{Rapha\"{e}l Carroy}

\address{
  Rapha\"{e}l Carroy \\
  Kurt G\"{o}del Research Center for Mathematical Logic \\
  Universit\"{a}t Wien \\
  W\"{a}hringer Stra{\ss}e 25 \\
  1090 Wien \\
  Austria
}

\email{raphael.carroy@univie.ac.at}

\urladdr{
  http://www.logique.jussieu.fr/~carroy/indexeng.html
}

\author[B.D. Miller]{Benjamin D. Miller}

\address{
  Benjamin D. Miller \\
  Kurt G\"{o}del Research Center for Mathematical Logic \\
  Universit\"{a}t Wien \\
  W\"{a}hringer Stra{\ss}e 25 \\
  1090 Wien \\
  Austria
 }

\email{benjamin.miller@univie.ac.at}

\urladdr{
  http://www.logic.univie.ac.at/benjamin.miller
}

\thanks{The authors were supported in part by FWF Grants
  P28153 and P29999.}
  


\title[The first Baire class]{Bases for functions beyond the first Baire
  class}

\maketitle

\section*{Introduction}

A topological space is \definedterm{analytic} if it is a continuous
image of a closed subset of $\Bairespace$. A subset of a topological
space is \definedterm{\Borel} if it is in the $\sigma$-algebra
generated by open sets, \definedterm{\Fsigma} if it is a union of
countably-many closed sets, and \definedterm{\Gdelta} if it is an
intersection of countably-many open sets.

Suppose that $X$ and $Y$ are topological spaces. Given a family
$\Gamma$ of subsets of $X$, a function $\phi \from X \to Y$ is
\definedterm{$\Gamma$-measurable} if $\preimage{\phi}{V} \in
\Gamma$ for every open set $V \subseteq Y$. A function is
\definedterm{\Borel} if it is \Borel-measurable, \definedterm{\Baire
class one} if it is \Fsigma-measurable, and \definedterm
{$\sigma$-continuous with closed witnesses} if its domain is the
union of countably-many closed sets on which it is
continuous. A result of \Jayne-\Rogers (see \cite[Theorem 1]
{JayneRogers}) ensures that a function from an analytic metric
space to a separable metric space has this property if and only if it
is \Gdelta-measurable.

A \definedterm{quasi-order} on a set $Z$ is a reflexive transitive
binary relation $\le$ on $Z$. A set $B \subseteq Z$ is a \definedterm
{basis} under $\le$ for $Z$ if $\forall z \in Z \exists b \in B \ b \le z$.

A \definedterm{closed continuous embedding} of $\phi \from X \to Y$
into $\phi' \from X' \to Y'$ consists of a pair of closed continuous
embeddings $\pi_X \from X \to X'$ and $\pi_Y \from \closure{\image
{\phi}{X}} \to \closure{\image{\phi'}{X'}}$ such that $\phi' \composition
\pi_X = \pi_Y \composition \phi$. Note that the existence of such a pair
depends not only on the graphs of the functions $\phi$ and $\phi'$, but
on $Y$ as well, since different choices of $Y \supseteq \image{\phi}{X}$
can lead to different values of $\closure{\image{\phi}{X}}$. Here we
establish the following results.

\begin{introtheorem} \label{intro:Baireclassone}
  There is a twenty-four-element basis under closed continuous
  embeddability for the class of non-\Baire-class-one \Borel functions
  between analytic metric spaces.
\end{introtheorem}

\begin{introtheorem} \label{intro:sigmacontinuous}
  There is a twenty-seven-element basis under closed continuous
  embeddability for the class of
  non-$\sigma$-continuous-with-closed-wit\-nesses \Borel functions
  between analytic metric spaces.
\end{introtheorem}

In \S\ref{compactification}, we discuss the compactification
$\compactifiedextendedBairespace$ of $\extendedBairespace$
underlying our arguments, as well as the corresponding
compactification $\compactifiedBairespace$ of $\Bairespace$. In \S\ref
{meetembeddings}, we discuss the endomorphisms of $\Bairetree$
underlying our arguments. In \S\ref{Bairespace}, we provide a
three-element basis for the class of \Baire measurable functions
from $\Bairespace$ to separable metric spaces. In \S\ref
{sigmacontinuous}, we provide a three-element basis for the class
of non-$\sigma$-continuous-with-closed-witnesses \Baire-class-one
functions from analytic metric spaces to separable metric spaces.
In \S\ref{Bairetree}, we provide an eight-element basis for the class
of all functions from $\compactifiedBairespace \setminus \Bairespace$
to analytic metric spaces. And in \S\ref{Baireclassone}, we establish
Theorems \ref{intro:Baireclassone} and \ref{intro:sigmacontinuous}.

\section{A compactification of $\extendedBairespace$}
  \label{compactification}

We use $s \concatenation t$ to denote the \definedterm{concatenation}
of sequences $s$ and $t$, and we say that $s$ is an \definedterm
{initial segment} of $t$, or $s \extendedby t$, if there exists $s'$ for
which $t = s \concatenation s'$. Endow the set
$\compactifiedextendedBairespace = \extendedBairespace \union \set
{t \concatenation \sequence{\infty}}[t \in \Bairetree]$ with the smallest
topology with respect to which the sets of the form $\set{t}$ and
$\extensions{t} = \set{c \in \compactifiedextendedBairespace}[t
\extendedby c]$, where $t \in \Bairetree$, are clopen.

\begin{proposition} \label{compactification:basis}
  The family $\calB$ of sets of the form $\set{t}$ and $\extensions{t}
  \setminus (\set{t} \union \union[j < i][\extensions{t \concatenation
  \sequence{j}}])$, where $i \in \N$ and $t \in \Bairetree$, is a clopen
  basis for $\compactifiedextendedBairespace$.
\end{proposition}

\begin{propositionproof}
  Let $\tau$ be the topology generated by $\calB$. As every set in
  $\calB$ is clearly clopen, it is sufficient to show that the sets $\set{t}$
  and $\extensions{t}$ are $\tau$-clopen for all $t \in \Bairetree$. As
  these sets are clearly $\tau$-open, we need only show that they are
  $\tau$-closed. As $\extensions{t \concatenation \sequence{i}}$ is
  $\tau$-closed in $\extensions{t}$ for all $i \in \N$ and $t \in
  \Bairetree$, a straightforward induction shows that $\extensions{t}$
  is $\tau$-closed for all $t \in \Bairetree$. As $\set{t}$ is $\tau$-closed
  in $\extensions{t}$ for all $t \in \Bairetree$, it follows that $\set{t}$ is
  $\tau$-closed for all $t \in \Bairetree$.
\end{propositionproof}

\begin{proposition} \label{compactification:compact}
  The space $\compactifiedextendedBairespace$ is compact.
\end{proposition}

\begin{propositionproof}
  Suppose, towards a contradiction, that there is an open cover
  $\calU$ of $\compactifiedextendedBairespace$ with no finite
  subcover.
  
  \begin{lemma} \label{compactification:compact:extension}
    Suppose that $t \in \Bairetree$ and no finite set $\calV \subseteq
    \calU$ covers $\extensions{t}$. Then there exists $j \in \N$ such
    that no finite set $\calV \subseteq \calU$ covers $\extensions{t
    \concatenation \sequence{j}}$.
  \end{lemma}
  
  \begin{lemmaproof}
    Fix $U \in \calU$ containing $t \concatenation \sequence{\infty}$.
    Proposition \ref{compactification:basis} then yields $i \in \N$ with
    $\extensions{t} \setminus (\set{t} \union \union[j < i][\extensions{t
    \concatenation \sequence{j}}]) \subseteq U$, in which case no
    finite set $\calV \subseteq \calU$ covers $\union[j < i][\extensions
    {t \concatenation \sequence{j}}]$, and it follows that there exists
    $j < i$ for which no finite set $\calV \subseteq \calU$ covers
    $\extensions{t \concatenation \sequence{j}}$.
  \end{lemmaproof}
  
  By recursively applying Lemma \ref
  {compactification:compact:extension}, we obtain $b \in
  \Bairespace$ such that for no $i \in \N$ is there a finite set
  $\calV \subseteq \calU$ covering $\extensions{\restriction{b}
  {i}}$. But Proposition \ref{compactification:basis} implies that
  every open neighborhood of $b$ contains some $\extensions
  {\restriction{b}{i}}$.
\end{propositionproof}

Given a countable set $I$ and a topological space $X$, we say
that a sequence $\sequence{x_i}[i \in I] \in \functions{I}{X}$
\definedterm{converges} to a point $x \in X$, or $x_i \goesto x$,
if for every open neighborhood $U$ of $x$ there are only finitely
many $i \in I$ with $x_i \notin U$. When $I$ and $X$ are
equipped with partial orders $\le_I$ and $\le_X$, we say that
$\sequence{x_i}[i \in I]$ is \definedterm{decreasing} if $i \le_I j
\implies x_j \le_X x_i$ for all $i, j \in I$.

\begin{proposition}
  The space $\compactifiedextendedBairespace$ has a
  compatible ultrametric.
\end{proposition}

\begin{propositionproof}
  Fix a decreasing sequence $\sequence{\epsilon_t}[t \in \Bairetree]$
  of positive real numbers converging to zero. Set $d(a, a) = 0$ for
  all $a \in \compactifiedextendedBairespace$, as well as $d(a, b) =
  \max \set{\epsilon_t}[t \in \set{\restriction{a}{\min(\length{a}, i(a, b))},
  \restriction{b}{\min(\length{b}, i(a, b))}} \intersection \Bairetree]$ for
  all distinct $a, b \in \compactifiedextendedBairespace$, where $i(a,
  b) = \min \set{i \in \N}[\restriction{a}{i} \neq \restriction{b}{i}]$.
  
  To see that $d$ is an ultrametric, suppose that $a, b, c \in
  \compactifiedextendedBairespace$ are pairwise distinct. Observe
  that if $i(a, c) < \max \set{i(a, b), i(b, c)}$, then $d(a, c) \in \set{d(b,
  c), d(a, b)}$, so $d(a, c) \le \max \set{d(a, b), d(b, c)}$. And if $i(a,
  c) = \max \set{i(a, b), i(b, c)}$, then setting $i = i(a, b) = i(a, c) =
  i(b, c)$, it follows that
  \begin{align*}
    d(a, c)
      & = \max \set{\epsilon_t}[t \in \set{\restriction{a}{i}, \restriction{c}
        {i}} \intersection \Bairetree] \\
      & \le \max \set{\epsilon_t}[t \in \set{\restriction{a}{i}, \restriction
        {b}{i}, \restriction{c}{i}} \intersection \Bairetree] \\
      & = \max \set{d(a, b),
        d(b, c)}.
  \end{align*}
  And if $i(a, c) > \max \set{i(a, b), i(b, c)}$, then setting $\epsilon =
  d(a, b) = d(b, c)$ and $t = \restriction{a}{i(a, b)} = \restriction{c}{i(b,
  c)}$, it follows that $d(a, c) \le \epsilon_t \le \epsilon$, and therefore
  $d(a, c) \le \max \set{d(a, b), d(b, c)}$.
  
  As $\set{t} = \ball{t}{\epsilon_t}$ and $\extensions{t} \setminus \set{t}
  = \ball{\extensions{t} \setminus \set{t}}{\epsilon_t}$ for all $t \in
  \Bairetree$, and $\extensions{t} \setminus (\set{t} \union \union[j \le i]
  [\extensions{t \concatenation \sequence{j}}]) = \ball{\extensions{t}
  \setminus (\set{t} \union \union[j \le i][\extensions{t \concatenation
  \sequence{j}}])}{\min(\set{\epsilon_{t \concatenation \sequence{j}}}[j
  \le i])}$ for all $i \in \N$ and $t \in \Bairetree$, Proposition \ref
  {compactification:basis} ensures that every open subset of
  $\compactifiedextendedBairespace$ is $d$-open.
  
  Given $b \in \Bairespace$ and $\epsilon > 0$, fix $i \in \N$ with
  $\epsilon_{\restriction{b}{i}} < \epsilon$, set $t = \restriction{b}{i}$,
  and note that $\extensions{t} \subseteq \ball{b}{\epsilon}$. Given
  $t \in \Bairetree$ and $\epsilon > 0$, fix $i \in \N$ with $\epsilon_{t
  \concatenation \sequence{j}} < \epsilon$ for all $j \ge i$, and
  observe that $\extensions{t} \setminus (\set{t} \union \union[j < i]
  [\extensions{t \concatenation \sequence{j}}]) \subseteq \ball{t
  \concatenation \sequence{\infty}}{\epsilon}$. Thus every $d$-open
  subset of $\compactifiedextendedBairespace$ is open.
\end{propositionproof}

It follows that $\compactifiedextendedBairespace$ is \Polish. As the 
space $\compactifiedBairespace = \compactifiedextendedBairespace
\setminus \Bairetree$ is a perfect subset of
$\compactifiedextendedBairespace$, a result of \Brouwer's ensures
that it is homeomorphic to $\Cantorspace$ (see, for example, \cite
[Theorem 7.4]{Kechris}).

\section{Meet embeddings} \label{meetembeddings}

The \definedterm{meet} of sequences $s, t \in \Bairetree$ is the
sequence $r = s \meet t$ of maximal length for which $r \extendedby
s$ and $r \extendedby t$. A \definedterm{$\meet$-embedding} is an
injection $\pi \from \Bairetree \to \Bairetree$ such that $\pi(s \meet t)
= \pi(s) \meet \pi(t)$ for all $s, t \in \Bairetree$.

\begin{proposition} \label{meetembeddings:equivalent}
  Suppose that $\pi \from \Bairetree \to \Bairetree$. Then $\pi$ is a
  $\meet$-embedding if and only if the following conditions hold:
  \begin{enumerate}
    \item $\forall i \in \N \forall t \in \Bairetree \ \pi(t) \strictlyextendedby
      \pi(t \concatenation \sequence{i})$.
    \item $\forall i, j \in \N \forall t \in \Bairetree$ \\
      \hspace*{5pt} $(i \neq j \implies \pi(t \concatenation \sequence{i})
        (\length{\pi(t)}) \neq \pi(t \concatenation \sequence{j})(\length
          {\pi(t)}))$.
  \end{enumerate}
\end{proposition}

\begin{propositionproof}
  Suppose first that $\pi$ is a $\meet$-embedding. To see that
  condition (1) holds, observe that if $i \in \N$ and $t \in \Bairetree$,
  then $\pi(t) = \pi(t) \meet \pi(t \concatenation \sequence{i})$, so
  $\pi(t) \extendedby \pi(t \concatenation \sequence{i})$, thus $\pi(t)
  \strictlyextendedby \pi(t \concatenation \sequence{i})$. And to see
  that condition (2) holds, observe that if $i, j \in \N$ are distinct and $t
  \in \Bairetree$, then $\pi(t) = \pi(t \concatenation \sequence{i})
  \meet \pi(t \concatenation \sequence{j})$, so $\pi(t \concatenation
  \sequence{i})(\length{\pi(t)}) \neq \pi(t \concatenation \sequence{j})
  (\length{\pi(t)})$.
  
  Suppose now that $\pi$ satisfies conditions (1) and (2). To see that
  $\pi$ is a $\meet$-embedding, suppose that $s, t \in \Bairetree$
  are distinct, and define $r = s \meet t$. By reversing the roles of
  $s$ and $t$ if necessary, we can assume that $\length{s} > \length
  {r}$, so $\pi(r) \strictlyextendedby \pi(s)$, thus either $r = t$ or
  ($\length{t} > \length{r}$ and $\pi(s)(\length{\pi(r)}) \neq \pi(t)
  (\length{\pi(r)})$). In both cases, it follows that $\pi(s) \neq
  \pi(t)$ and $\pi(r) = \pi(s) \meet \pi(t)$.
\end{propositionproof}

\begin{remark}
  In particular, it follows that if $\pi \from \Bairetree \to \Bairetree$ has
  the property that $\pi(t) \concatenation \sequence{i} \extendedby
  \pi(t \concatenation \sequence{i})$ for all $i \in \N$ and $t \in \Bairetree$,
  then $\pi$ is a $\meet$-embedding.
\end{remark}

There is a simple but useful means of amalgamating appropriately
indexed families of $\meet$-embeddings.

\begin{proposition} \label{meetembeddings:composition}
  Suppose that $\sequence{\pi_t}[t \in \Bairetree]$ is a sequence of
  $\meet$-embeddings with the property that $\image{\pi_t}
  {\Bairetree} \subseteq \extensions{t}$ for all $t \in \Bairetree$.
  Then the function $\pi \from \Bairetree \to \Bairetree$ given by
  $\pi(t) = (\product[n \le \length{t}][\pi_{\restriction{t}{n}}])(t)$ is also
  a $\meet$-embedding.
\end{proposition}

\begin{propositionproof}
  Note that if $i \in \N$ and $t \in \Bairetree$, then $t \concatenation
  \sequence{i} \extendedby \pi_{t \concatenation \sequence{i}}(t
  \concatenation \sequence{i})$, so Proposition \ref
  {meetembeddings:equivalent} ensures that $(\product[n \le \length
  {t}][\pi_{\restriction{t}{n}}])(t \concatenation \sequence{i})
  \extendedby \pi(t \concatenation \sequence{i})$, thus $\pi(t)
  \strictlyextendedby (\product[n \le \length{t}][\pi_{\restriction{t}{n}}])
  (t \concatenation \sequence{i}) \extendedby \pi(t \concatenation
  \sequence{i})$. It also implies that if $i \neq j$, then $(\product[n
  \le \length{t}][\pi_{\restriction{t}{n}}])(t \concatenation \sequence{i})
  (\length{\pi(t)}) \neq (\product[n \le \length{t}][\pi_{\restriction{t}
  {n}}])(t \concatenation \sequence{j})(\length{\pi(t)})$, so $\pi(t
  \concatenation \sequence{i})(\length{\pi(t)}) \neq \pi(t
  \concatenation \sequence{j})(\length{\pi(t)})$. One last application
  of Proposition \ref{meetembeddings:equivalent} therefore ensures
  that $\pi$ is a $\meet$-embedding.
\end{propositionproof}

We next consider the connection between $\meet$-embeddings
and closed continuous embeddings.

\begin{proposition}
  Every $\meet$-embedding $\pi \from \Bairetree \to \Bairetree$ has
  a unique extension to a (necessarily injective) continuous map
  $\extension{\pi} \from \compactifiedextendedBairespace \to
  \compactifiedextendedBairespace$, given by $\extension{\pi}(b) =
  \union[i \in \N][\pi(\restriction{b}{i})]$ and $\extension{\pi}(t
  \concatenation \sequence{\infty}) = \pi(t) \concatenation \sequence
  {\infty}$ for all $b \in \Bairespace$ and $t \in \Bairetree$.
\end{proposition}

\begin{propositionproof}
  Suppose that $\extension{\pi} \from
  \compactifiedextendedBairespace \to
  \compactifiedextendedBairespace$ is a continuous extension of
  $\pi$. If $b \in \Bairespace$, then $\restriction{b}{i} \goesto b$, and
  since $\sequence{\pi(\restriction{b}{i})}[i \in \N]$ is strictly
  increasing by Proposition \ref{meetembeddings:equivalent}, it
  follows that $\extension{\pi}(b) = \union[i \in \N][\pi(\restriction{b}
  {i})]$. If $t \in \Bairetree$, then $t \concatenation \sequence{i}
  \goesto t \concatenation \sequence{\infty}$, and since $\pi(t) =
  \pi(t \concatenation \sequence{i}) \meet \pi(t \concatenation
  \sequence{j})$ for all distinct $i, j \in \N$, it follows that $\extension
  {\pi}(t \concatenation \sequence{\infty}) = \pi(t) \concatenation
  \sequence{\infty}$.
  
  To see that these constraints actually define a continuous function,
  note that if $t \in \Bairetree$, then either $\preimage{\extension
  {\pi}}{\extensions{t}} = \emptyset$ or there exists $s \in \Bairetree$
  of minimal length with $t \extendedby \pi(s)$, in which case
  $\preimage{\extension{\pi}}{\extensions{t}} = \extensions{s}$.
  
  To see that $\extension{\pi}$ is injective, it is enough to check that
  its restriction to $\Bairespace$ is injective. Towards this end,
  suppose that $a, b \in \Bairespace$ are distinct, fix $i \in \N$ least
  for which $a(i) \neq b(i)$, set $t = \restriction{a}{i} = \restriction{b}
  {i}$, and observe that $\pi(t \concatenation \sequence{a(i)})(\length
  {\pi(t)}) \neq \pi(t \concatenation \sequence{b(i)})(\length{\pi(t)})$
  by Proposition \ref{meetembeddings:equivalent}, thus $\extension
  {\pi}(a)$ and $\extension{\pi}(b)$ are distinct.
\end{propositionproof}

\begin{remark} \label{meetembeddings:composition:remark}
  It follows that the extension associated with the composition of two
  $\meet$-embeddings is the composition of their extensions.
\end{remark}

Given a function $\phi \from X \to Y$ and sets $X' \subseteq X$ and
$Y' \supseteq \image{\phi}{X'}$, let $\restriction{\phi}{X'}[Y']$ denote
the function $\psi \from X' \to Y'$ given by $\phi(x) = \psi(x)$ for all $x
\in X'$. Compactness ensures that if $\pi$ is a $\meet$-embedding,
then $\extension{\pi}$ and $\restriction{\extension{\pi}}
{\compactifiedBairespace}$ are closed continuous embeddings. The
following observations show that so too are $\restriction{\extension
{\pi}}{\Bairespace}[\Bairespace]$ and $\restriction{\extension{\pi}}
{\compactifiedBairespace \setminus \Bairespace}
[\compactifiedBairespace \setminus \Bairespace]$.

\begin{proposition} \label{meetembeddings:closed}
  Suppose that $\pi \from \Bairetree \to \Bairetree$ is a
  $\meet$-embedding. Then $\restriction{\extension{\pi}}
  {\Bairespace}[\Bairespace]$ is closed.
\end{proposition}

\begin{propositionproof}
  It is sufficient to show that every sequence $\sequence{b_n}[n \in
  \N]$ of elements of $\Bairespace$ for which $\sequence{\extension
  {\pi}(b_n)}[n \in \N]$ converges to an element of $\Bairespace$ is
  itself convergent to an element of $\Bairespace$. As $\sequence
  {\restriction{\extension{\pi}(b_n)}{i}}[n \in \N]$ is eventually constant
  for all $i \in \N$, a simple induction shows that $\sequence
  {\restriction{b_n}{i}}[n \in \N]$ is also eventually constant for all $i \in
  \N$, so $\sequence{b_n}[n \in \N]$ converges to an element of
  $\Bairespace$.
\end{propositionproof}

\begin{proposition} \label{meetembeddings:closed:difference}
  Suppose that $\pi \from \Bairetree \to \Bairetree$ is a
  $\meet$-embedding. Then $\restriction{\extension{\pi}}
  {\compactifiedBairespace \setminus \Bairespace}
  [\compactifiedBairespace \setminus \Bairespace]$ is closed.
\end{proposition}

\begin{propositionproof}
  It is sufficient to show that every sequence $\sequence{s_n}[n \in
  \N]$ of elements of $\Bairetree$ such that $\sequence{\pi(s_n)}[n
  \in \N]$ converges to $t \concatenation \sequence{\infty}$ for some
  $t \in \Bairetree$ has a subsequence converging to an element of
  $\compactifiedBairespace \setminus \Bairespace$. By passing to a
  subsequence, we can assume that $\pi(s_m) \meet \pi(s_n) = t$ for
  all distinct $m, n \in \N$. Let $s$ be the $\extendedby$-minimal
  element of $\Bairetree$ for which $t \extendedby \pi(s)$. Then
  $s_m \meet s_n = s$ for all distinct $m, n  \in \N$, thus $s_n
  \goesto s \concatenation \sequence{\infty}$.
\end{propositionproof}

A set $T \subseteq \Bairetree$ is \definedterm
{$\extendedby$-dense} if $\forall s \in \Bairetree \exists t \in T \ s
\extendedby t$. More generally, a set $T \subseteq \Bairetree$ is
\definedterm{$\extendedby$-dense} below $r \in \Bairetree$ if
$\forall s \in \Bairetree \exists t \in T \ r \concatenation s
\extendedby t$.

\begin{proposition} \label{meetembeddings:ramsey}
  Suppose that $T \subseteq \Bairetree$. Then there is a
  $\meet$-embedding $\pi \from \Bairetree \to \Bairetree$ such that
  $\image{\pi}{\Bairetree} \subseteq T$ or $\image{\pi}{\Bairetree}
  \subseteq \setcomplement{T}$. 
\end{proposition}

\begin{propositionproof}
  Fix $S \in \set{T, \setcomplement{T}}$ which is $\extendedby$-dense
  below some $s \in \Bairetree$, and recursively
  construct a function $\pi \from \Bairetree \to \extensions{s}
  \intersection S$ with the property that $\pi(t) \concatenation \sequence{i}
  \extendedby \pi(t \concatenation \sequence{i})$ for all $i \in \N$
  and $t \in \Bairetree$.
\end{propositionproof}

\begin{proposition} \label{meetembeddings:category}
  Suppose that $C \subseteq \Bairespace$ is a non-meager set with
  the \Baire property. Then there is a $\meet$-embedding $\pi \from
  \Bairetree \to \Bairetree$ with the property that $\image{\extension
  {\pi}}{\Bairespace} \subseteq C$.
\end{proposition}

\begin{propositionproof}
  Fix $s \in \Bairetree$ for which $C$ is comeager in $\extensions
  {s} \intersection \Bairespace$, as well as dense open sets $U_n
  \subseteq \extensions{s} \intersection \Bairespace$ with the
  property that $\intersection[n \in \N][U_n] \subseteq C$. Set $T_n
  = \set{t \in \Bairetree}[\extensions{t} \intersection \Bairespace
  \subseteq U_n]$ for all $n \in \N$, and recursively construct a
  function $\pi \from \Bairetree \to \extensions{s} \intersection
  \Bairetree$ such that $\image{\pi}{\Bairespace[n]} \subseteq
  T_n$ for all $n \in \N$ and $\pi(t) \concatenation \sequence{i}
  \extendedby \pi(t \concatenation \sequence{i})$ for all $i \in \N$
  and $t \in \Bairetree$.
\end{propositionproof}

\section{Baire measurable functions on $\Bairespace$}
\label{Bairespace}

Here we provide a basis for the class of \Baire measurable functions
from $\Bairespace$ to separable metric spaces.

\begin{proposition} \label{Bairespace:continuous}
  Suppose that $X$ is a second countable topological space and
  $\phi \from \Bairespace \to X$ is \Baire measurable. Then there is
  a $\meet$-embedding $\pi \from \Bairetree \to \Bairetree$ for
  which $\phi \composition \extension{\pi}$ is continuous. 
\end{proposition}

\begin{propositionproof}
  Fix a comeager set $C \subseteq \Bairespace$ on which $\phi$
  is continuous, and appeal to Proposition \ref
  {meetembeddings:category} to obtain a $\meet$-embedding $\pi
  \from \Bairetree \to \Bairetree$ with the property that $\image
  {\extension{\pi}}{\Bairespace} \subseteq C$.
\end{propositionproof}

\begin{proposition} \label{Bairespace:diameter}
  Suppose that $X$ is a metric space and $\phi \from \Bairespace \to
  X$ is continuous. Then there is a $\meet$-embedding $\pi \from
  \Bairetree \to \Bairetree$ with the property that $\diameter{\image
  {\phi}{\extensions{\pi(t)}}} \goesto 0$.
\end{proposition}

\begin{propositionproof}
  Fix a sequence $\sequence{\epsilon_t}[t \in \Bairetree]$ of positive
  real numbers converging to zero, note that the continuity of $\phi$
  ensures that for all $t \in \Bairetree$ the set $T_t = \set{s \in
  \Bairetree}[\diameter{\image{\phi}{\extensions{s}}} < \epsilon_t]$ is
  $\extendedby$-dense, and recursively construct a function $\pi
  \from \Bairetree \to \Bairetree$ such that $\pi(t) \in T_t$ for all $t \in
  \Bairetree$ and $\pi(t) \concatenation \sequence{i} \extendedby \pi(t
  \concatenation \sequence{i})$ for all $i \in \N$ and $t \in \Bairetree$.
\end{propositionproof}

Given a countable set $I$ and a topological space $X$, we say
that a sequence $\sequence{X_i}[i \in I]$ of subsets of $X$
\definedterm{converges} to a point $x \in X$, or $X_i \goesto x$, if
for every open neighborhood $U$ of $x$, all but finitely many $i \in I$
have the property that $X_i \subseteq U$. We say that $\sequence
{X_i}[i \in I]$ is \definedterm{discrete} if for all $x \in X$ there is an
open neighborhood $U$ of $x$ such that all but finitely many $i \in I$
have the property that $U \intersection X_i = \emptyset$.

\begin{proposition} \label{Bairespace:atmostone}
  Suppose that $X$ is a metric space and $\phi \from \Bairespace \to
  X$ has the property that $\diameter{\image{\phi}{\extensions{t
  \concatenation \sequence{i}}}} \goesto 0$ for all $t \in \Bairetree$.
  Then there is a $\meet$-embedding $\pi \from \Bairetree \to
  \Bairetree$ such that $\sequence{\image{\phi}{\extensions{\pi(t
  \concatenation \sequence{i})}}}[i \in \N]$ is convergent or discrete
  for all $t \in \Bairetree$.
\end{proposition}

\begin{propositionproof}
  For each $t \in \Bairetree$, the fact that $\diameter{\image{\phi}
  {\extensions{t \concatenation \sequence{i}}}} \goesto 0$ ensures
  that there is an injection $\iota_t \from \N \to \N$ for which
  $\sequence{\image{\phi}{\extensions{t \concatenation \sequence
  {\iota_t(i)}}}}[i \in \N]$ is convergent or discrete. Define $\pi \from
  \Bairetree \to \Bairetree$ by choosing $\pi(\emptysequence) \in
  \Bairetree$ arbitrarily and setting $\pi(t \concatenation \sequence
  {i}) = \pi(t) \concatenation \sequence{\iota_{\pi(t)}(i)}$ for all $i \in
  \N$ and $t \in \Bairetree$.
\end{propositionproof}

We say that a function $\phi \from X \to Y$ is \definedterm{nowhere
constant} if there is no non-empty open set $U \subseteq X$ on which
$\phi$ is constant.

\begin{proposition} \label{Bairespace:disjointness}
  Suppose that $X$ is a metric space and $\phi \from \Bairespace \to
  X$ is continuous and nowhere constant. Then there is a
  $\meet$-embedding $\pi \from \Bairetree \to \Bairetree$ such that
  \begin{equation*}
    \textstyle
    \forall i \in \N \forall t \in \Bairetree \ \closure{\image{\phi}
      {\extensions{\pi(t \concatenation \sequence{i})}}} \intersection
        \closure{\union[j \in \N \setminus \set{i}][\image{\phi}
          {\extensions{\pi(t \concatenation \sequence{j})}}]} = \emptyset.
  \end{equation*}
\end{proposition}

\begin{propositionproof}
  Clearly each $\image{\phi}{\extensions{t}}$ is infinite.
  
  \begin{lemma}
    For all $t \in \Bairetree$, there is a function $\iota_t \from \N \to
    \Bairetree \setminus \set{\emptysequence}$ such that $\sequence
    {\iota_t(i)(0)}[i \in \N]$ is injective and the closures of $\image
    {\phi}{\extensions{t \concatenation \iota_t(i)}}$ and
    $\union[j \in \N \setminus \set{i}][\image{\phi}{\extensions{t
    \concatenation \iota_t(j)}}]$ are disjoint for all $i \in \N$.
  \end{lemma}
  
  \begin{lemmaproof}
    As each $\image{\phi}{\extensions{t \concatenation \sequence{i}}}$ is
    infinite, there are extensions $b_i \in \Bairespace$
    of $t \concatenation \sequence{i}$ such that $\phi(b_i)
    \notin \set{\phi(b_j)}[j < i]$ for all $i \in \N$. Fix a subsequence
    $\sequence{a_i}[i \in \N]$ of $\sequence{b_i}[i \in \N]$ for which
    $\set{\phi(a_i)}[i \in \N]$ is discrete. For each $i \in \N$, fix
    $\epsilon_i > 0$ such that $\phi(a_j) \notin \ball{\phi(a_i)}
    {\epsilon_i}$ for all $j \in \N \setminus \set{i}$, as well as $\iota_t(i)
    \in \Bairetree \setminus \set{\emptysequence}$ with $t
    \concatenation \iota_t(i) \extendedby a_i$ and $\image{\phi}
    {\extensions{t \concatenation \iota_t(i)}} \subseteq \ball{\phi(a_i)}
    {\epsilon_i / 3}$.

    Suppose, towards a contradiction, that there exists $i \in \N$ for
    which some $x \in X$ is in the closures of $\image{\phi}
    {\extensions{t \concatenation \iota_t(i)}}$ and $\union
    [j \in \N \setminus \set{i}][\image{\phi}{\extensions{t \concatenation
    \iota_t(j)}}]$. Then there exist $j \in \N \setminus \set
    {i}$ and $y \in \image{\phi}{\extensions{t \concatenation
    \iota_t(j)}}$ with the property that $d(x, y) \le \epsilon_i / 3$, in
    which case
    \begin{align*}
      d(\phi(a_i), \phi(a_j))
        & \le d(\phi(a_i), x) + d(x, y) + d(y, \phi(a_j)) \\
        & < \epsilon_i / 3 + \epsilon_i / 3 + \epsilon_j / 3 \\
        & \le \max \set{\epsilon_i, \epsilon_j},
    \end{align*}
    so $\phi(a_i) \in \ball{\phi(a_j)}{\epsilon_j}$ or $\phi(a_j) \in \ball
    {\phi(a_i)}{\epsilon_i}$, a contradiction.
  \end{lemmaproof}
  
  Define $\pi \from \Bairetree \to \Bairetree$ by choosing $\pi
  (\emptysequence) \in \Bairetree$ arbitrarily and setting $\pi(t
  \concatenation \sequence{i}) = \pi(t) \concatenation \iota_{\pi(t)}(i)$
  for all $i \in \N$ and $t \in \Bairetree$.
\end{propositionproof}

We now obtain our main result stabilizing the topological behavior
of \Baire measurable functions from $\Bairespace$ to separable
metric spaces.

\begin{theorem} \label{Bairespace:closedcontinuousembedding}
  Suppose that $X$ is a separable metric space and $\phi \from
  \Bairespace \to X$ is \Baire measurable. Then there is a
  $\meet$-embedding $\pi \from \Bairetree \to \Bairetree$ such that
  $\phi \composition \extension{\pi}$ is constant or extends to a
  closed continuous embedding on $\Bairespace$ or
  $\compactifiedBairespace$.
\end{theorem}

\begin{theoremproof}
  By Remark \ref{meetembeddings:composition:remark}, we are
  free to replace $\phi$ by its composition with the extension of any
  $\meet$-embedding. For example, by Proposition \ref
  {Bairespace:continuous}, we can assume that $\phi$ is
  continuous.
  
  If there exists $s \in \Bairetree$ for which $\restriction{\phi}
  {\extensions{s}}$ is constant, then define $\pi \from \Bairetree
  \to \Bairetree$ by $\pi(t) = s \concatenation t$ for all $t \in
  \Bairetree$, so $\phi \composition \extension{\pi}$ is constant.
  Otherwise, Propositions \ref{meetembeddings:ramsey}, \ref
  {Bairespace:diameter}, \ref{Bairespace:atmostone}, and \ref
  {Bairespace:disjointness} yield a $\meet$-embedding $\pi \from
  \Bairetree \to \Bairetree$ such that $\diameter{\image{\phi}
  {\extensions{\pi(t)}}} \goesto 0$, $\sequence{\image{\phi}
  {\extensions{\pi(t \concatenation \sequence{i})}}}[i \in \N]$ is
  convergent for all $t \in \Bairetree$ or discrete for all $t \in
  \Bairetree$, and
  \begin{equation*}
    \textstyle
    \forall i \in \N \forall t \in \Bairetree \ \closure{\image{\phi}
      {\extensions{\pi(t \concatenation \sequence{i})}}} \intersection
        \closure{\union[j \in \N  \setminus \set{i}][\image{\phi}
          {\extensions{\pi(t \concatenation \sequence{j})}}]} = \emptyset.
  \end{equation*}
  As $\image{\extension{\pi}}{\extensions{t}} \subseteq \extensions
  {\pi(t)}$ for all $t \in \Bairetree$, it follows that
  \begin{equation*}
    \textstyle
    \forall i \in \N \forall t \in \Bairetree \ \closure{\image{(\phi
      \composition \extension{\pi})}{\extensions{t \concatenation
        \sequence{i}}}} \intersection \closure{\union[j \in \N \setminus
          \set{i}][\image{(\phi \composition \extension{\pi})}{\extensions{t
            \concatenation \sequence{j}}}]} = \emptyset.
  \end{equation*}
  So by replacing $\phi$ with $\phi \composition \extension{\pi}$, we
  can assume that $\diameter{\image{\phi}{\extensions{t}}} \goesto
  0$, $\sequence{\image{\phi}{\extensions{t \concatenation
  \sequence{i}}}}[i \in \N]$ is convergent for all $t \in \Bairetree$ or
  discrete for all $t \in \Bairetree$, and
  \begin{equation} \tag{$\dagger$}
    \textstyle
    \forall i \in \N \forall t \in \Bairetree \ \closure{\image{\phi}
      {\extensions{t \concatenation \sequence{i}}}} \intersection
        \closure{\union[j \in \N  \setminus \set{i}][\image{\phi}
          {\extensions{t \concatenation \sequence{j}}}]} = \emptyset.
  \end{equation}
  
  To see that $\phi$ is injective, note that if $a, b \in \Bairespace$
  are distinct, then there is a least $i \in \N$ for which $a(i) \neq b(i)$.
  Setting $t = \restriction{a}{i} = \restriction{b}{i}$, it follows from
  $(\dagger)$ that $\image{\phi}{\extensions{t \concatenation
  \sequence{a(i)}}}$ and $\image{\phi}{\extensions{t \concatenation
  \sequence{b(i)}}}$ are disjoint, thus $\phi(a)$ and $\phi(b)$ are
  distinct.

  We next check that if $\sequence{\image{\phi}{\extensions{t
  \concatenation \sequence{i}}}}[i \in \N]$ is discrete for all $t
  \in \Bairetree$, then $\phi$ is a closed continuous embedding. It is
  sufficient to show that every sequence $\sequence{b_n}[n \in \N]$
  of elements of $\Bairespace$ for which $\sequence{\phi(b_n)}[n
  \in \N]$ converges to some $x \in X$ is itself convergent. But a
  straightforward recursive argument yields $b \in \Bairespace$ such
  that $x$ is in the closure of $\image{\phi}{\extensions{\restriction{b}
  {i}}}$ for all $i \in \N$, so $(\dagger)$ ensures that $x$ is not in the
  closure of $\union[j \in \N \setminus \set{b(i)}][\image{\phi}
  {\extensions{\restriction{b}{i} \concatenation \sequence{j}}}]$ for
  all $i \in \N$, thus $\sequence{\restriction{b_n}{i}}[n \in \N]$ is
  eventually constant with value $\restriction{b}{i}$ for all $i \in \N$,
  hence $b_n \goesto b$.
 
  It remains to check that if $\sequence{\image{\phi}{\extensions
  {t \concatenation \sequence{i}}}}[i \in \N]$ is convergent for all $t
  \in \Bairetree$, then the extension of $\phi$ to
  $\compactifiedBairespace$ given by $\extension{\phi}(t
  \concatenation \sequence{\infty}) = \lim_{i \goesto \infty} \image
  {\phi}{\extensions{t \concatenation \sequence{i}}}$ for all $t \in
  \Bairetree$ is a closed continuous embedding. To
  see that $\extension{\phi}$ is injective, note that if $c, d \in
  \compactifiedBairespace$ are distinct, then there is a least $i \in
  \N$ with $c(i) \neq d(i)$. By reversing the roles of $c$ and $d$ if
  necessary, we can assume that $c(i) \neq \infty$. Set $t =
  \restriction{c}{i} = \restriction{d}{i}$, and appeal to $(\dagger)$ to
  see that $\extension{\phi}(c)$ is in the closure of $\image{\phi}
  {\extensions{t \concatenation \sequence{c(i)}}}$ but $\extension
  {\phi}(d)$ is not, so $\extension{\phi}(c) \neq \extension{\phi}(d)$.
  To see that $\extension{\phi}$ is continuous, suppose that $c \in
  \compactifiedBairespace$ and $U$ is an open neighborhood of
  $\extension{\phi}(c)$, and fix an open neighborhood $V$ of
  $\extension{\phi}(c)$ whose closure is contained in $U$. If $c \in
  \Bairespace$, then there exists $i \in \N$ for which $\image{\phi}
  {\extensions{\restriction{c}{i}}} \subseteq V$, thus $\extensions
  {\restriction{c}{i}}$ is an open neighborhood of $c$ whose image
  under $\extension{\phi}$ is contained in $U$. Otherwise, there
  exists $t \in \Bairetree$ for which $c = t \concatenation \sequence
  {\infty}$, as well as $i \in \N$ for which $\image{\phi}{\extensions{t}
  \setminus \union[j < i][\extensions{t \concatenation \sequence{j}}]}
  \subseteq V$. Then $\extensions{t} \setminus \union[j < i]
  [\extensions{t \concatenation \sequence{j}}]$ is an open
  neighborhood of $c$ whose image under $\extension{\phi}$ is
  contained in $U$.
\end{theoremproof}

For each topological space $X$, let $\constant{X}$ denote the
unique function from $X$ to the trivial topological space $\set{\infty}$.
Given topological spaces $X \subseteq Y$, define $\inclusion{X}{Y}
\from X \to Y$ by $\inclusion{X}{Y}(x) = x$ for all $x \in X$.

\begin{proposition} \label{Bairespace:basis}
  Suppose that $X$ is a separable metric space, $\phi \from
  \Bairespace \to X$ is \Baire measurable, $\pi \from \Bairetree \to
  \Bairetree$ is a $\meet$-embedding, and $\phi \composition
  \extension{\pi}$ is constant or extends to a closed continuous
  embedding on $\Bairespace$ or $\compactifiedBairespace$. Then
  there exist $\phi_0 \in \set{\constant{\Bairespace}} \union \set
  {\inclusion{\Bairespace}{Z}}[Z \in \set{\Bairespace,
  \compactifiedBairespace}]$ and $\psi \from \closure{\image{\phi_0}
  {\Bairespace}} \to \closure{\image{\phi}{\Bairespace}}$ with the
  property that $\pair{\restriction{\extension{\pi}}{\Bairespace}
  [\Bairespace]}{\psi}$ is a closed continuous embedding of $\phi_0$
  into $\phi$.
\end{proposition}

\begin{propositionproof}
  If $\phi \composition \extension{\pi}$ is constant, then set $\phi_0
  = \constant{\Bairespace}$ and let $\psi$ be the unique function
  from $\image{\constant{\Bairespace}}{\Bairespace}$ to $\image
  {(\phi \composition \extension{\pi})}{\Bairespace}$. If $\phi
  \composition \extension{\pi}$ extends to a closed continuous
  embedding $\psi$ on $Z \in \set{\Bairespace,
  \compactifiedBairespace}$, then set $\phi_0 = \inclusion
  {\Bairespace}{Z}$.
\end{propositionproof}

\section{\Baire-class-one functions that are not
  $\sigma$-continuous with closed witnesses}
  \label{sigmacontinuous}

Here we strengthen \cite
[Theorem 3.1]{Solecki} by providing a basis for the class of
non-$\sigma$-continuous-with-closed-witnesses \Baire-class-one
functions from analytic metric spaces to separable metric spaces.

\begin{proposition} \label{sigmacontinuous:limit}
  Suppose that $X$ is a metric space and $\phi \from
  \compactifiedBairespace \to X$ has the property that $\restriction
  {\phi}{\Bairespace}$ is continuous. Then there is a
  $\meet$-embedding $\pi \from \Bairetree \to \Bairetree$ such that
  either $\closure{\image{(\phi \composition \extension{\pi})}
  {\Bairespace}} \intersection \closure{\image{(\phi \composition
  \extension{\pi})}{\compactifiedBairespace \setminus \Bairespace}}
  = \emptyset$ or $\phi \composition \extension{\pi}$ is continuous
  at every point of $\Bairespace$.
\end{proposition}

\begin{propositionproof}
  We can assume that there is no $s \in \Bairetree$ with the property
  that $\inf \set{d(\phi(s \concatenation b), \phi(s \concatenation t
  \concatenation \sequence{\infty}))}[b \in \Bairespace \mathand t \in
  \Bairetree] > 0$, since otherwise the $\meet$-embedding $\pi \from
  \Bairetree \to \Bairetree$ given by $\pi(t) = s \concatenation t$ for
  all $t \in \Bairetree$ has the property that $\closure{\image{(\phi
  \composition \extension{\pi})}{\Bairespace}} \intersection \closure
  {\image{(\phi \composition \extension{\pi})}{\compactifiedBairespace
  \setminus \Bairespace}} = \emptyset$.
  
  \begin{lemma}
    Suppose that $\epsilon > 0$ and $s \in \Bairetree$. Then there
    exists $t \in \Bairetree$ with $d(\phi(s \concatenation t
    \concatenation b), \phi(s \concatenation t \concatenation
    \sequence{\infty})) < \epsilon$ for all $b \in \Bairespace$.
  \end{lemma}
  
  \begin{lemmaproof}
    Fix $\delta < \epsilon$ and $u \in \Bairetree$ with $\diameter
    {\image{\phi}{\extensions{s \concatenation u} \intersection
    \Bairespace}} < \delta$, and $b \in \Bairespace$ and $v \in
    \Bairetree$ with $d(\phi(s \concatenation u \concatenation
    b), \phi(s \concatenation u \concatenation v \concatenation
    \sequence{\infty})) < \epsilon - \delta$, and set $t = u
    \concatenation v$.
  \end{lemmaproof}
  
  Fix a sequence $\sequence{\epsilon_n}[n \in \N]$ of positive real
  numbers converging to zero, and recursively construct a function
  $\pi \from \Bairetree \to \Bairetree$ with the property that $d
  (\phi(\pi(t) \concatenation b), \phi(\pi(t) \concatenation \sequence
  {\infty})) < \epsilon_{\length{t}}$ for all $b \in \Bairespace$ and $t
  \in \Bairetree$, and $\pi(t) \concatenation \sequence{i} \extendedby
  \pi(t \concatenation \sequence{i})$ for all $i \in \N$ and $t \in
  \Bairetree$.
\end{propositionproof}

We say that a metric space is \definedterm{$\epsilon$-discrete} if all
distinct points have distance at least $\epsilon$ from one another.

\begin{proposition} \label{sigmacontinuous:closeddiscrete:epsilon}
  Suppose that $X$ is a metric space, $\phi \from
  \compactifiedBairespace \setminus \Bairespace \to X$,
  $\epsilon > 0$, and $t \in \Bairetree$. Then there is a
  $\meet$-embedding $\pi \from \Bairetree \to \extensions{t}
  \intersection \Bairetree$ with the property that $\phi \composition
  \extension{\pi}$ is an injection into an $\epsilon$-discrete set or
  $\image{(\phi \composition \extension{\pi})}
  {\compactifiedBairespace \setminus \Bairespace}$ is contained in
  the $\epsilon$-ball around a point of $\image{\phi}{\extensions{t}}$.
\end{proposition}

\begin{propositionproof}
  If for no finite set $F \subseteq \image{\phi}
  {\compactifiedBairespace \setminus \Bairespace}$ and extension
  $u$ of $t$ is it the case that $\image{\phi}{\extensions{u}}
  \subseteq \ball{F}{\epsilon}$, then fix an enumeration $\sequence
  {t_n}[n \in \N]$ of $\Bairetree$ with the property that $t_m
  \extendedby t_n \implies m \le n$ for all $m, n \in \N$, and
  recursively construct $\pi \from \Bairetree \to \extensions{t}
  \intersection \Bairetree$ such that $\phi(\pi(t_n) \concatenation
  \sequence{\infty}) \notin \ball{\set{\phi(\pi(t_m) \concatenation
  \sequence{\infty})}[m < n]}{\epsilon}$ and $\pi(t_n') \concatenation
  \sequence{n} \extendedby \pi(t_n)$ for all $n > 0$, where $t_n'$ is
  the maximal proper initial segment of $t_n$.

  Otherwise, there exists $x \in \image{\phi}{\compactifiedBairespace
  \setminus \Bairespace}$ with the property that the set $S = \set{s
  \in \Bairetree}[\phi(s \concatenation \sequence{\infty}) \in \ball{x}
  {\epsilon}]$ is $\extendedby$-dense below some extension $u$ of
  $t$, in which case we can recursively construct a function $\pi \from
  \Bairetree \to \extensions{u} \intersection S$ with the property that
  $\pi(v) \concatenation \sequence{i} \extendedby \pi(v \concatenation
  \sequence{i})$ for all $i \in \N$ and $v \in \Bairetree$.
\end{propositionproof}

\begin{proposition} \label{sigmacontinuous:closeddiscrete}
  Suppose that $X$ is a metric space and $\phi \from
  \compactifiedBairespace \setminus \Bairespace \to X$. Then there
  is a $\meet$-embedding $\pi \from \Bairetree \to \Bairetree$ such
  that $\phi \composition \extension{\pi}$ is an injection into an
  $\epsilon$-discrete set for some $\epsilon > 0$ or $\diameter
  {\image{(\phi \composition \extension{\pi})}{\extensions{t}}}
  \goesto 0$.
\end{proposition}

\begin{propositionproof}
  Suppose that for no $\epsilon > 0$ is there a $\meet$-embedding
  $\pi \from \Bairetree \to \Bairetree$ such that $\phi \composition
  \extension{\pi}$ is an injection into an $\epsilon$-discrete set, fix a
  sequence $\sequence{\epsilon_t}[t \in \Bairetree]$ of positive real
  numbers converging to zero, and recursively apply Proposition \ref
  {sigmacontinuous:closeddiscrete:epsilon} to the functions $\phi_t =
  \phi \composition \product[n < \length{t}][\extension{\pi_{\restriction
  {t}{n}}}]$ to obtain $\meet$-embeddings $\pi_t \from \Bairetree \to
  \extensions{t} \intersection \Bairetree$ such that $\image{(\phi
  \composition \product[n \le \length{t}][\extension{\pi_{\restriction{t}
  {n}}}])}{\compactifiedBairespace \setminus \Bairespace}$ is
  contained in an $\epsilon_t$-ball for all $t \in \Bairetree$. Let $\pi$
  be the $\meet$-embedding obtained from applying Proposition \ref
  {meetembeddings:composition} to $\sequence{\pi_t}[t \in
  \Bairetree]$, and observe that $\diameter{\image{(\phi
  \composition \extension{\pi})}{\extensions{t}}} \goesto 0$.
\end{propositionproof}

Define $\discrete \from \compactifiedBairespace \setminus
\Bairespace \to \Bairetree$ by setting $\discrete(t \concatenation
\sequence{\infty}) = t$ for all $t \in \Bairetree$. Let
$\onepointcompactification{\Bairetree} = \Bairetree \union \set
{\infty}$ denote the \definedterm{one-point compactification} of
$\Bairetree$.

\begin{theorem} \label{sigmacontinuous:basis}
  Suppose that $X$ is an analytic metric space, $Y$ is a separable
  metric space, and $\phi \from X \to Y$ is a \Baire-class-one function
  that is not $\sigma$-continuous with closed witnesses. Then there
  exists $\phi_0 \in \set{\constant{\Bairespace}} \union \set{\inclusion
  {\Bairespace}{Z}}[Z \in \set{\Bairespace, \compactifiedBairespace}]$
  for which there is a closed continuous embedding of $\phi_0 \union
  \discrete$ into $\phi$.
\end{theorem}

\begin{theoremproof}
  By the \Jayne-\Rogers theorem (see, for example, \cite[Theorem
  1]{JayneRogers}), we can assume that $\phi$ is not
  \Gdelta-measurable. \Hurewicz's dichotomy theorem for \Fsigma
  sets then yields a closed continuous embedding $\psi \from
  \compactifiedBairespace \to X$ with $\closure{\image{(\phi
  \composition \psi)}{\Bairespace}} \intersection \image{(\phi
  \composition \psi)}{\compactifiedBairespace \setminus \Bairespace}
  = \emptyset$ (see, for example, \cite[Theorem 4.2]{CarroyMillerSoukup}).
  As $\pair{\psi}{\identity{\closure{\image{(\phi \composition \psi)}
  {\compactifiedBairespace}}}}$ is a closed continuous embedding
  of $\phi \composition \psi$ into $\phi$, by replacing the latter with
  the former, we can assume that $X = \compactifiedBairespace$
  and $\closure{\image{\phi}{\Bairespace}} \intersection \image{\phi}
  {\compactifiedBairespace \setminus \Bairespace} = \emptyset$.
  
  By Proposition \ref{Bairespace:continuous}, there is a
  $\meet$-embedding $\pi \from \Bairetree \to \Bairetree$ for which
  $\restriction{(\phi \composition \extension{\pi})}{\Bairespace}$ is
  continuous. By composing $\pi$ with the $\meet$-embedding
  given by Proposition \ref{sigmacontinuous:limit}, we can assume
  that $\closure{\image{(\phi \composition \extension{\pi})}
  {\Bairespace}} \intersection \closure{\image{(\phi \composition
  \extension{\pi})}{\compactifiedBairespace \setminus \Bairespace}}
  = \emptyset$ or $\phi \composition \extension{\pi}$ is continuous
  at every point of $\Bairespace$. As $\phi$ is \Baire class one, the
  former possibility would imply that the pre-images of $\closure{\image
  {(\phi \composition \extension{\pi})}{\Bairespace}}$ and $\closure
  {\image{(\phi \composition \extension{\pi})}{\compactifiedBairespace
  \setminus \Bairespace}}$ under $\phi \composition \extension{\pi}$
  are disjoint dense \Gdelta subsets of $\compactifiedBairespace$, so
  the latter holds. By Proposition \ref{sigmacontinuous:closeddiscrete},
  we can assume that either there exists $\epsilon > 0$ for which
  $\restriction{(\phi \composition \extension{\pi})}
  {\compactifiedBairespace \setminus \Bairespace}$ is an injection into
  an $\epsilon$-discrete set, or $\diameter{\image{(\phi \composition
  \extension{\pi})}{\extensions{t} \intersection (\compactifiedBairespace
  \setminus \Bairespace)}} \goesto 0$. As the former possibility
  contradicts the facts that $\image{(\phi \composition \extension{\pi})}
  {\Bairespace} \intersection \image{(\phi \composition \extension{\pi})}
  {\compactifiedBairespace \setminus \Bairespace} = \emptyset$ and
  $\image{(\phi \composition \extension{\pi})}{\Bairespace} \subseteq
  \closure{\image{(\phi \composition \extension{\pi})}
  {\compactifiedBairespace \setminus \Bairespace}}$, it follows that
  the latter holds. By applying Proposition \ref
  {sigmacontinuous:closeddiscrete:epsilon} with any $\epsilon > 0$
  and $t \in \Bairetree$, but replacing the given metric on $X$ by one
  with respect to which all pairs of distinct points have distance at
  least $\epsilon$ from one another, we can assume that $\restriction
  {(\phi \composition \extension{\pi})}{\compactifiedBairespace
  \setminus \Bairespace}$ is either constant or injective.
  
  \begin{lemma} \label{sigmacontinuous:basis:sequences}
    Suppose that $\sequence{s_n}[n \in \N]$ is an injective sequence
    of elements of $\Bairetree$ and $\sequence{b_n}[n \in \N]$ is a
    sequence of elements of $\Bairespace$ such that $s_n
    \extendedby b_n$ for all $n \in \N$. Then $\metric{X}((\phi \composition
    \extension{\pi})(b_n), (\phi \composition \extension{\pi})(s_n
    \concatenation \sequence{\infty})) \goesto 0$.
  \end{lemma}
  
  \begin{lemmaproof}
    Simply note that $(\phi \composition \extension{\pi})(b_n) \in
    \closure{\image{(\phi \composition \extension{\pi})}{\extensions
    {s_n} \intersection (\compactifiedBairespace \setminus
    \Bairespace)}}$ for all $n \in \N$ and $\diameter{\image{(\phi
    \composition \extension{\pi})}{\extensions{s_n} \intersection
    (\compactifiedBairespace \setminus \Bairespace)}} \goesto 0$.
  \end{lemmaproof}
  
  Along with the facts that $\image{(\phi
  \composition \extension{\pi})}{\Bairespace} \intersection \image
  {(\phi \composition \extension{\pi})}{\compactifiedBairespace
  \setminus \Bairespace} = \emptyset$ and $\image{(\phi
  \composition \extension{\pi})}{\Bairespace} \subseteq \closure
  {\image{(\phi \composition \extension{\pi})}
  {\compactifiedBairespace \setminus \Bairespace}}$, Lemma \ref
  {sigmacontinuous:basis:sequences} ensures that $\restriction{(\phi
  \composition \extension{\pi})}{\compactifiedBairespace \setminus
  \Bairespace}$ is not constant, and is therefore injective. Along with
  the fact that $\closure{\image{(\phi \composition \extension{\pi})}
  {\Bairespace}} \intersection \image{(\phi \composition \extension
  {\pi})}{\compactifiedBairespace \setminus \Bairespace} = \emptyset$,
  Lemma \ref{sigmacontinuous:basis:sequences} ensures that
  $\image{(\phi \composition \extension{\pi})}{\compactifiedBairespace
  \setminus \Bairespace}$ is discrete.
  
  By Theorem \ref{Bairespace:closedcontinuousembedding}, we can
  assume that $\restriction{(\phi \composition \extension{\pi})}
  {\Bairespace}$ is constant or extends to a closed continuous
  embedding on $\Bairespace$ or $\compactifiedBairespace$.
  
  We will now complete the proof by showing that there exist
  $\phi_0 \in \set{\constant{\Bairespace}} \union \set{\inclusion
  {\Bairespace}{Z}}[Z \in \set{\Bairespace, \compactifiedBairespace}]$
  and $\psi \from \closure{\image{\phi_0}{\compactifiedBairespace}}
  \union \Bairetree \to \closure{\image{\phi}{X}}$ for which $\pair
  {\restriction{\extension{\pi}}{\compactifiedBairespace}
  [\compactifiedBairespace]}{\psi}$ is a closed continuous embedding
  of $\phi_0 \union \discrete$ into $\phi$.
    
  If $\restriction{(\phi \composition \extension{\pi})}{\Bairespace}$
  is constant with value $y \in Y$, then set $\phi_0 = \constant
  {\Bairespace}$, and note that the extension $\psi$ of $\phi
  \composition \extension{\pi} \composition \inverse{\discrete}$ to
  $\onepointcompactification{\Bairetree}$ given by $\psi(\infty) = y$ is
  injective. As Lemma \ref{sigmacontinuous:basis:sequences} ensures
  that $(\phi \composition \extension{\pi})(s_n \concatenation \sequence
  {\infty}) \goesto y$ for every injective sequence $\sequence{s_n}
  [n \in \N]$ of elements of $\Bairetree$, it follows that $\psi$ is
  continuous, so the compactness of $\onepointcompactification
  {\Bairetree}$ ensures that $\psi$ is a closed continuous embedding.
  
  If $\restriction{(\phi \composition \extension{\pi})}{\Bairespace}$ is
  a closed continuous embedding, then set $\phi_0 = \inclusion
  {\Bairespace}{\Bairespace}$, and note that the extension $\psi$ of
  $\phi \composition \extension{\pi} \composition \inverse{\discrete}$
  to $\extendedBairespace$ given by $\restriction{\psi}{\Bairespace}
  = \restriction{(\phi \composition \extension{\pi})}{\Bairespace}$ is a
  continuous injection. To see that it is closed, it is enough to show that
  every injective sequence $\sequence{a_n}[n \in \N]$ of elements of
  $\extendedBairespace$ for which $\sequence{\psi(a_n)}[n \in \N]$
  converges to some point $y \in Y$ has a subsequence converging
  to a point of $\Bairespace$. As $\compactifiedextendedBairespace$
  is compact, by passing to a subsequence, we can assume that
  $\sequence{a_n}[n \in \N]$ converges to a point of
  $\compactifiedextendedBairespace$. As every point of $\Bairetree$
  is isolated, it therefore converges to a point of
  $\compactifiedBairespace$. And if there exists $t \in \Bairetree$ for
  which $a_n \goesto t \concatenation \sequence{\infty}$, then there
  are extensions $b_n \in \Bairespace$ of $a_n$ for all $n \in \N$, in
  which case $b_n \goesto t \concatenation \sequence{\infty}$ and
  $\psi(b_n) \goesto y$ by Lemma \ref
  {sigmacontinuous:basis:sequences}. Fix $n \in \N$ sufficiently large
  that $(\phi \composition \extension{\pi})(b_m) \neq y$ for all $m \ge
  n$, and observe that $\set{b_m}[m \ge n]$ is a closed subset of
  $\Bairespace$ whose image under $\phi \composition \extension
  {\pi}$ is not closed, contradicting the fact that $\restriction{(\phi
  \composition \extension{\pi})}{\Bairespace}$ is closed.
  
  If $\restriction{(\phi \composition \extension{\pi})}{\Bairespace}$
  extends to a closed continuous embedding $\psi'$ on
  $\compactifiedBairespace$, then set $\phi_0 = \inclusion
  {\Bairespace}{\compactifiedBairespace}$, and note that the
  extension $\psi$ of $\phi \composition \extension{\pi}
  \composition \inverse{\discrete}$ to
  $\compactifiedextendedBairespace$ given by $\restriction{\psi}
  {\compactifiedBairespace} = \restriction{\psi'}
  {\compactifiedBairespace}$ is injective. To see that it is continuous,
  suppose that $\sequence{t_n}[n \in \N]$ is an injective sequence of
  elements of $\Bairetree$ converging to $t \concatenation \sequence
  {\infty}$ for some $t \in \Bairetree$, fix $b_n \in \extensions{t_n}
  \intersection \Bairespace$ for all $n \in \N$, and observe that the
  continuity of $\psi'$ ensures that $\psi(b_n) \goesto \psi(t
  \concatenation \sequence{\infty})$, thus Lemma \ref
  {sigmacontinuous:basis:sequences} implies that $\psi(t_n) \goesto
  \psi(t \concatenation \sequence{\infty})$. As
  $\compactifiedextendedBairespace$ is compact, it follows that
  $\psi$ is a closed continuous embedding.
\end{theoremproof}

\section{Functions on $\compactifiedBairespace \setminus
  \Bairespace$} \label{Bairetree}
  
Here we provide a basis for the class of all functions
from $\compactifiedBairespace \setminus \Bairespace$ to analytic
metric spaces.

\begin{proposition} \label{Bairetree:pointavoidance}
  Suppose that $X$ is a topological space, $\phi \from
  \compactifiedBairespace \setminus \Bairespace \to X$ is injective,
  and $x \in X$. Then there is a $\meet$-embedding $\pi \from
  \Bairetree \to \Bairetree$ such that $x \notin \image{(\phi
  \composition \extension{\pi})}{\compactifiedBairespace \setminus
  \Bairespace}$.
\end{proposition}

\begin{propositionproof}
  Fix $s \in \Bairetree$ such that $x \notin \image{\phi}{\extensions
  {s}}$, and define $\pi \from \Bairetree \to \Bairetree$ by $\pi(t) =
  s \concatenation t$ for all $t \in \Bairetree$.
\end{propositionproof}

\begin{proposition} \label{Bairetree:convergentordiscrete}
  Suppose that $X$ is a metric space and $\phi \from
  \compactifiedBairespace \setminus \Bairespace \to X$. Then there
  is a $\meet$-embedding $\pi \from \Bairetree \to \Bairetree$ with
  the property that $\sequence{(\phi \composition \extension{\pi})(t
  \concatenation \sequence{i, \infty})}[i \in \N]$ is convergent or
  $\set{(\phi \composition \extension{\pi})(t \concatenation \sequence
  {i, \infty})}[i \in \N]$ is closed and discrete for all $t \in \Bairetree$.
\end{proposition}

\begin{propositionproof}
  For each $t \in \Bairetree$, there is an injection $\iota_t \from \N \to
  \N$ for which $\sequence{\phi(t \concatenation \sequence{\iota_t(i),
  \infty})}[i \in \N]$ is convergent or $\set{\phi(t \concatenation
  \sequence{\iota_t(i), \infty})}[i \in \N]$ is closed and discrete. Define
  $\pi \from \Bairetree \to \Bairetree$ by choosing $\pi
  (\emptysequence) \in \Bairetree$ arbitrarily and setting $\pi(t
  \concatenation \sequence{i}) = \pi(t) \concatenation \sequence
  {\iota_{\pi(t)}(i)}$ for all $i \in \N$ and $t \in \Bairetree$, and note
  that $(\phi \composition \extension{\pi})(t \concatenation
  \sequence{i, \infty}) = \phi(\pi(t \concatenation \sequence{i})
  \concatenation \sequence{\infty}) = \phi(\pi(t) \concatenation
  \sequence{\iota_{\pi(t)}(i), \infty})$ for all $i \in \N$ and $t \in
  \Bairespace$.
\end{propositionproof}

\begin{proposition} \label{Bairetree:discrete:F}
  Suppose that $X$ is a metric space, $\phi \from
  \compactifiedBairespace \setminus \Bairespace \to X$, $F
  \subseteq X$ is finite, and $t \in \Bairetree$. Then there is a
  $\meet$-embedding $\pi \from \Bairetree \to \extensions{t}
  \intersection \Bairetree$ such that either $\sequence{(\phi
  \composition \extension{\pi})(u \concatenation \sequence{\infty})}[u
  \in \Bairetree]$ converges to an element of $F$ or the closure of
  $\image{(\phi \composition \extension{\pi})}
  {\compactifiedBairespace \setminus \Bairespace}$ is disjoint from
  $F$.
\end{proposition}

\begin{propositionproof}
  If the set $S_\epsilon = \set{s \in \Bairetree}[\phi(s \concatenation
  \sequence{\infty}) \in \ball{F}{\epsilon}]$ is $\extendedby$-dense
  below $t$ for all $\epsilon > 0$, then there exist an extension $u$
  of $t$ and $x \in F$ such that the set $S_{\epsilon, x} = \set{s \in
  \Bairetree}[\phi(s \concatenation \sequence{\infty}) \in \ball{x}
  {\epsilon}]$ is $\extendedby$-dense below $u$ for all $\epsilon >
  0$. Fix a sequence $\sequence{\epsilon_v}[v \in \Bairetree]$ of
  positive real numbers converging to zero, and recursively
  construct a function $\pi \from \Bairetree \to \extensions{u}
  \intersection \Bairetree$ such that $\pi(v) \in S_{\epsilon_v, x}$
  for all $v \in \Bairetree$ and $\pi(v) \concatenation \sequence{i}
  \extendedby \pi(v \concatenation \sequence{i})$ for all $i \in \N$
  and $v \in \Bairetree$, and observe that $(\phi \composition
  \extension{\pi})(v \concatenation \sequence{\infty}) \goesto x$.  
      
  Otherwise, fix $\epsilon > 0$ and an extension $u$ of $t$ with
  the property that $\extensions{u} \intersection S_\epsilon =
  \emptyset$, define $\pi \from \Bairetree \to \extensions{u}
  \intersection \Bairetree$ by $\pi(v) = u \concatenation v$, and
  note that the closure of $\image{(\phi \composition \extension
  {\pi})}{\compactifiedBairespace \setminus \Bairespace}$ is
  disjoint from $F$.
\end{propositionproof}

For the rest of this section, it will be convenient to fix an enumeration
$\sequence{t_n}[n \in \N]$ of $\Bairetree$ such that $t_m
\extendedby t_n \implies m \le n$ for all $m, n \in \N$.

\begin{proposition} \label{Bairetree:discrete}
  Suppose that $X$ is a metric space and $\phi \from
  \compactifiedBairespace \setminus \Bairespace \to X$. Then there
  is a $\meet$-embedding $\pi \from \Bairetree \to \Bairetree$ with
  the property that $\sequence{(\phi \composition \extension{\pi})(t
  \concatenation \sequence{\infty})}[t \in \Bairetree]$ converges or for
  no natural numbers $m < n$ is $(\phi \composition \extension{\pi})
  (t_m \concatenation \sequence{\infty})$ or a limit point of $\set{(\phi
  \composition \extension{\pi})(t_m \concatenation \sequence{i,
  \infty})}[i \in \N]$ in the closure of $\image{(\phi \composition
  \extension{\pi})}{\extensions{t_n}}$.
\end{proposition}

\begin{propositionproof}
  Suppose that for no $\meet$-embedding $\pi \from \Bairetree
  \to \Bairetree$ is the sequence $\sequence{(\phi \composition
  \extension{\pi})(t \concatenation \sequence{\infty})}[t \in
  \Bairetree]$ convergent. By Proposition \ref
  {Bairetree:convergentordiscrete}, we can assume that $\sequence
  {\phi(t \concatenation \sequence{i, \infty})}[i \in \N]$ is convergent
  or $\set{\phi(t \concatenation \sequence{i, \infty})}[i \in \N]$ is
  closed and discrete for all $t \in \Bairetree$. By recursively
  applying Lemma \ref{Bairetree:discrete:F} to the functions
  $\phi_t = \phi \composition \product[k < \length{t}][\extension
  {\pi_{\restriction{t}{k}}}]$, we obtain $\meet$-embeddings $\pi_t
  \from \Bairetree \to \extensions{t} \intersection \Bairetree$ such
  that for no natural numbers $m < n$ is $(\phi \composition \prod_{k \le \length{t_m}}
  \extension{\pi_{\restriction{t_m}{k}}})(t_m \concatenation \sequence
  {\infty})$ or a limit point of $\set{(\phi \composition \prod_{k \le
  \length{t_m}} \extension{\pi_{\restriction{t_m}{k}}})(t_m
  \concatenation \sequence{i, \infty})}[i \in \N]$ in the closure of
  $\image{(\phi \composition \prod_{k \le \length{t_n}} \extension
  {\pi_{\restriction{t_n}{k}}})}{\extensions{t_n}}$. Let $\pi$ be the
  $\meet$-embedding obtained from applying Proposition \ref
  {meetembeddings:composition} to $\sequence{\pi_t}[t \in
  \Bairetree]$, and observe that for no natural numbers $m < n$ is it
  the case that $(\phi \composition \extension{\pi})(t_m
  \concatenation \sequence{\infty})$ or a limit point of $\set{(\phi
  \composition \extension{\pi})(t_m \concatenation \sequence{i,
  \infty})}[i \in \N]$ in the closure of $\image{(\phi \composition
  \extension{\pi})}{\extensions{t_n}}$.
\end{propositionproof}

\begin{theorem}
  \label{Bairetree:closedcontinuousembedding}
  Suppose that $X$ is an analytic metric space and $\phi \from
  \compactifiedBairespace \setminus \Bairespace \to X$. Then there
  is a $\meet$-embedding $\pi \from \Bairetree \to \Bairetree$ such
  that $\phi \composition \extension{\pi}$ is constant, $\phi
  \composition \extension{\pi}$ extends to a closed continuous
  embedding on $\compactifiedBairespace \setminus \Bairespace$
  or $\compactifiedBairespace$, or $\phi \composition \extension
  {\pi} \composition \inverse{\discrete}$ extends to a closed
  continuous embedding on $\Bairetree$, $\onepointcompactification
  {\Bairetree}$, $\compactifiedextendedBairespace \setminus \Bairespace$, $\extendedBairespace$, or
  $\compactifiedextendedBairespace$.
\end{theorem}
  
\begin{theoremproof}
  As before, we will repeatedly precompose $\phi$ with appropriate
  $\meet$-embeddings, albeit this time so as to stabilize the behavior of the function 
$\psi = \phi \composition \inverse{\discrete}$, as opposed to that of the function $\phi$ itself. By applying
  Proposition \ref{sigmacontinuous:closeddiscrete:epsilon} with any
  $\epsilon > 0$ and $t \in \Bairetree$, but replacing the given metric
  on $X$ by one with respect to which all pairs of distinct points have
  distance at least $\epsilon$ from one another, we can assume that
  $\psi$ is either constant or injective. As $\phi$ is constant in the former
  case, we can assume that we are in the latter.

  By Proposition \ref{sigmacontinuous:closeddiscrete}, we can
  ensure that $\image{\psi}{\Bairetree}$ is closed and discrete or
  $\diameter{\image{\psi}{\extensions{t}}} \goesto 0$. As $\psi$ is a
  closed continuous embedding in the former case, we can assume
  that we are in the latter.
  
  Let $\extension{\psi}$ be the extension of $\psi$ to a partial
  function on $\compactifiedextendedBairespace$ given by
  $\extension{\psi}(b) = \lim_{i \goesto \infty} \psi(\restriction{b}{i})$
  and $\extension{\psi}(t \concatenation \sequence{\infty}) = \lim_{i
  \goesto \infty} \psi(t \concatenation \sequence{i})$ for all $b \in
  \Bairespace$ and $t \in \Bairetree$. By Proposition \ref
  {Bairetree:convergentordiscrete}, we can assume that $\set{\psi
  (t \concatenation \sequence{i})}[i \in \N]$ has a limit point $\implies
  t \concatenation \sequence{\infty} \in \domain{\extension{\psi}}$ for
  all $t \in \Bairetree$.
  
  As each point of $\Bairetree$ is isolated, $\diameter{\image{\psi}
  {\extensions{\restriction{b}{i}}}} \goesto 0$ for all $b \in
  \Bairespace$, and $\diameter{\image{\psi}{\extensions{t
  \concatenation \sequence{i}}}} \goesto 0$ for all $t \in \Bairetree$,
  it follows that $\extension{\psi}$ is continuous. To see that
  $\extension{\psi}$ is closed, it is sufficient show that every injective
  sequence $\sequence{c_n}[n \in \N]$ of points in the domain of
  $\extension{\psi}$ for which $\sequence{\extension{\psi}(c_n)}[n \in
  \N]$ is convergent has a subsequence converging to a point in the
  domain of $\extension{\psi}$. By passing to a subsequence, we can
  assume that the sequence converges to a point of
  $\compactifiedextendedBairespace$. As each point of $\Bairetree$
  is isolated, the sequence converges to a point of
  $\compactifiedBairespace$, so the facts that $\diameter{\image
  {\psi}{\extensions{\restriction{b}{i}}}} \goesto 0$ for all $b \in
  \Bairespace$, $\diameter{\image{\psi}{\extensions{t \concatenation
  \sequence{i}}}} \goesto 0$ for all $t \in \Bairetree$, and $\set{\psi(t
  \concatenation \sequence{i})}[i \in \N]$ has a limit point $\implies t
  \concatenation \sequence{\infty} \in \domain{\extension{\psi}}$ for
  all $t \in \Bairetree$ ensure that it converges to a point of the
  domain of $\extension{\psi}$.
  
  By Proposition \ref{meetembeddings:ramsey}, we can
  assume that one of the following holds:
  \begin{enumerate}
    \item $\compactifiedBairespace \setminus \Bairespace \subseteq
      \domain{\extension{\psi}} \mathand \forall t \in \Bairetree
        \ \extension{\psi}(t) = \extension{\psi}(t \concatenation
          \sequence{\infty})$.
    \item $\compactifiedBairespace \setminus \Bairespace \subseteq
      \domain{\extension{\psi}} \mathand \forall t \in \Bairetree
        \ \extension{\psi}(t) \neq \extension{\psi}(t \concatenation
          \sequence{\infty})$.
    \item $(\compactifiedBairespace \setminus \Bairespace)
      \intersection \domain{\extension{\psi}} = \emptyset$.
  \end{enumerate}
  As the domain of $\extension{\psi}$ is analytic, so too is its
  intersection with $\Bairespace$. It follows that the latter
  intersection has the \Baire property, so Proposition \ref
  {meetembeddings:category} allows us to assume that one of
  the following holds:
  \begin{enumerate}
    \renewcommand{\theenumi}{\alph{enumi}}
    \item The domain of $\extension{\psi}$ is disjoint from
      $\Bairespace$.
    \item The domain of $\extension{\psi}$ contains $\Bairespace$.
  \end{enumerate}
  In the special case that condition (b) holds, Theorem \ref
  {Bairespace:closedcontinuousembedding} allows us to assume
  that $\restriction{\extension{\psi}}{\Bairespace}$ is either constant
  or injective.
  
  Proposition \ref{Bairetree:discrete} allows us to assume that
  $\sequence{\psi(t)}[t \in \Bairetree]$ converges to some $x \in X$
  or for no natural numbers $m < n$ is $\psi(t_m)$ or $\extension
  {\psi}(t_m \concatenation \sequence{\infty})$ in the closure of
  $\image{\psi}{\extensions{t_n}}$. In the former case, Proposition
  \ref{Bairetree:pointavoidance} allows us to assume that $\image
  {\psi}{\Bairetree}$ is discrete, so the extension of $\psi$ to
  $\onepointcompactification{\Bairetree}$ sending $\infty$ to $x$ is a
  closed continuous embedding, thus we can assume that we are in
  the latter.
  
  \begin{lemma}
    Suppose that $c, d \in \domain{\extension{\psi}}$ are distinct
    but $\extension{\psi}(c) = \extension{\psi}(d)$. Then there exists $t
    \in \Bairetree$ such that $\set{c, d} = \set{t, t \concatenation
    \sequence{\infty}}$.
  \end{lemma}
  
  \begin{lemmaproof}
    To see that $\restriction{\extension{\psi}}{\compactifiedBairespace
    \setminus \Bairespace}$ is injective, observe that if $m < n$, both
    $t_m \concatenation \sequence{\infty}$ and $t_n \concatenation
    \sequence{\infty}$ are in the domain of $\extension{\psi}$, and
    moreover $\extension{\psi}(t_m \concatenation \sequence{\infty})
    = \extension{\psi}(t_n \concatenation \sequence{\infty})$, then
    $\extension{\psi}(t_m \concatenation \sequence{\infty})$ is in
    the closure of $\image{\psi}{\extensions{t_n}}$, a contradiction.
    
    To see that $\restriction{\extension{\psi}}{\Bairespace}$ is injective
    when $\Bairespace$ is contained in the domain of $\extension
    {\psi}$, note that otherwise it is constant, and let $x$ be this
    constant value. Then for each $t \in \Bairetree$, there is a
    sequence $\sequence{u_i}[i \in \N]$ of elements of $\Bairetree$
    such that $\psi(t \concatenation \sequence{i} \concatenation
    \sequence{u_i}) \goesto x$, so the fact that $\diameter{\image
    {\psi}{\extensions{t \concatenation \sequence{i}}}} \goesto 0$
    ensures that $\extension{\psi}(t \concatenation \sequence{\infty})
    = x$, contradicting the fact that $\restriction{\extension{\psi}}
    {\compactifiedBairespace \setminus \Bairespace}$ is injective.
    
    To see that $\image{\extension{\psi}}{\Bairespace} \intersection
    \image{\psi}{\Bairetree} = \emptyset$, note that if $b \in \domain
    {\extension{\psi}} \intersection \Bairespace$, $t \in \Bairetree$, and
    $\extension{\psi}(b) = \psi(t)$, then there exist $m < n$ with $t_m
    = t$ and $t_n \strictlyextendedby b$, so $\psi(t_m)$ is in the
    closure of $\image{\psi}{\extensions{t_n}}$, a contradiction.

    To see that $\image{\extension{\psi}}{\Bairespace} \intersection
    \image{\extension{\psi}}{\compactifiedBairespace \setminus
    \Bairespace} = \emptyset$, note that if $b \in \domain{\extension
    {\psi}} \intersection \Bairespace$, $t \in \Bairetree$, $t
    \concatenation \sequence{\infty} \in \domain{\extension{\psi}}$,
    and $\extension{\psi}(b) = \extension{\psi}(t \concatenation
    \sequence{\infty})$, then there exist $m < n$ with $t_m = t$ and
    $t_n \strictlyextendedby b$, in which case $\extension{\psi}(t_m
    \concatenation \sequence{\infty})$ is in the closure of $\image
    {\psi}{\extensions{t_n}}$, a contradiction.
    
    Observe finally that if $s, t \in \Bairetree$ are distinct, $t
    \concatenation \sequence{\infty} \in \domain{\extension{\psi}}$,
    and $\psi(s) = \extension{\psi}(t \concatenation \sequence
    {\infty})$, then there exist $m \neq n$ such that $t_m = s$ and
    $t_n = t$. Then $\psi(t_m)$ is in the closure of $\image{\psi}
    {\extensions{t_n}}$ and $\extension{\psi}(t_n \concatenation
    \sequence{\infty})$ is in $\image{\psi}{\extensions{t_m}}$, a
    contradiction.
  \end{lemmaproof}
  
  If (1a) or (1b) holds, then $\restriction{\extension{\psi}}
  {\compactifiedBairespace \setminus \Bairespace}$ or $\restriction
  {\extension{\psi}}{\compactifiedBairespace}$ is an extension of
  $\phi$ to a closed continuous embedding. If (2a), (2b), (3a), or (3b)
  holds, then $\extension{\psi}$ is an extension of $\psi$ to a closed
  continuous embedding on $\compactifiedextendedBairespace \setminus \Bairespace$,
  $\compactifiedextendedBairespace$, $\Bairetree$, or
  $\extendedBairespace$.
\end{theoremproof}

\begin{proposition} \label{Bairetree:basis}
  Suppose that $X$ is an analytic metric space, $\phi \from
  \compactifiedBairespace \setminus \Bairespace \to X$, $\pi \from
  \Bairetree \to \Bairetree$ is a $\meet$-embedding, and $\phi
  \composition \extension{\pi}$ is constant, $\phi \composition
  \extension{\pi}$ extends to a closed continuous embedding on
  $\compactifiedBairespace \setminus \Bairespace$ or
  $\compactifiedBairespace$, or $\phi \composition \extension{\pi}
  \composition \inverse{\discrete}$ extends to a closed continuous
  embedding on $\Bairetree$, $\onepointcompactification
  {\Bairetree}$, $\compactifiedextendedBairespace \setminus \Bairespace$, $\extendedBairespace$, or
  $\compactifiedextendedBairespace$. Then there exist $\phi_0 \in
  \set{\constant{\compactifiedBairespace \setminus \Bairespace}}
  \union \set{\inclusion{\compactifiedBairespace \setminus
  \Bairespace}{Z}}[Z \in \set{\compactifiedBairespace \setminus
  \Bairespace, \compactifiedBairespace}] \union \set{\inclusion
  {\Bairetree}{Z} \composition \discrete}[Z \in \set{\Bairetree,
  \onepointcompactification{\Bairetree}, \compactifiedextendedBairespace \setminus \Bairespace,
  \extendedBairespace, \compactifiedextendedBairespace}]$ and
  $\psi \from \closure{\image{\phi_0}{\compactifiedBairespace
  \setminus \Bairespace}} \to \closure{\image{\phi}
  {\compactifiedBairespace \setminus \Bairespace}}$ with the
  property that $\pair{\restriction{\extension{\pi}}
  {\compactifiedBairespace \setminus \Bairespace}
  [\compactifiedBairespace \setminus \Bairespace]}{\psi}$ is a closed
  continuous embedding of $\phi_0$ into $\phi$.
\end{proposition}

\begin{propositionproof}
  If $\phi \composition \extension{\pi}$ is constant, then set $\phi_0
  = \constant{\compactifiedBairespace \setminus \Bairespace}$ and
  let $\psi$ be the unique function from $\image{\constant
  {\compactifiedBairespace \setminus \Bairespace}}
  {\compactifiedBairespace \setminus \Bairespace}$ to $\image{(\phi
  \composition \extension{\pi})}{\compactifiedBairespace \setminus
  \Bairespace}$. If $\phi \composition \extension{\pi}$ extends to a
  closed continuous embedding $\psi$ on $Z \in \set
  {\compactifiedBairespace \setminus \Bairespace,
  \compactifiedBairespace}$, then set $\phi_0 = \inclusion
  {\compactifiedBairespace \setminus \Bairespace}{Z}$. And if $\phi
  \composition \extension{\pi} \composition \inverse{\discrete}$
  extends to a closed continuous embedding $\psi$ on $Z \in \set
  {\Bairetree, \onepointcompactification{\Bairetree}, \compactifiedextendedBairespace \setminus \Bairespace,
  \extendedBairespace, \compactifiedextendedBairespace}$, then
  set $\phi_0 = \inclusion{\Bairetree}{Z} \composition \discrete$.
\end{propositionproof}

\section{Borel functions that are not Baire class one}
  \label{Baireclassone}

Here we provide bases for the classes of non-\Baire-class-one \Borel
functions and non-$\sigma$-continuous-with-closed-witnesses \Borel
functions between analytic metric spaces.

\begin{proposition} \label{Baireclassone:disjoint}
  Suppose that $X$ is a metric space and $\phi \from
  \compactifiedBairespace \to X$ has the property that $\restriction
  {\phi}{\Bairespace}$ is continuous and $\image{\phi}{\Bairespace}
  \nsubseteq \closure{\image{\phi}{\compactifiedBairespace
  \setminus \Bairespace}}$. Then there is a $\meet$-embedding $\pi
  \from \Bairetree \to \Bairetree$ with the property that $\closure
  {\image{(\phi \composition \extension{\pi})}{\Bairespace}}
  \intersection \closure{\image{(\phi \composition \extension{\pi})}
  {\compactifiedBairespace \setminus \Bairespace}} = \emptyset$.
\end{proposition}

\begin{propositionproof}
    Fix  $b \in \Bairespace$ for which $\phi(b)$ is
    not in the closure of $\image{\phi}{\compactifiedBairespace
    \setminus \Bairespace}$. Then there is an open neighborhood $U$
    of $\phi(b)$ disjoint from $\image{\phi}{\compactifiedBairespace
    \setminus \Bairespace}$, as well as an open neighborhood $V$ of
    $\phi(b)$ whose closure is contained in $U$, in which case the
    continuity of $\restriction{\phi}{\Bairespace}$ yields a proper initial
    segment $s$ of $b$ for which $\image{\phi}{\extensions{s}
    \intersection \Bairespace} \subseteq V$. Then the
    $\meet$-embedding $\pi \from \Bairetree \to \Bairetree$ given by
    $\pi(t) = s \concatenation t$ for all $t \in \Bairetree$ is as desired.
\end{propositionproof}

Given $\phi_{\Bairespace} \from \Bairespace \to X$ and
$\phi_{\compactifiedBairespace \setminus \Bairespace} \from
\compactifiedBairespace \setminus \Bairespace \to Y$, let
$\phi_{\Bairespace} \disjointunion \phi_{\compactifiedBairespace
\setminus \Bairespace}$ denote the corresponding function from
$\compactifiedBairespace$ to the disjoint union $X \disjointunion Y$.

\begin{theorem} \label{Baireclassone:basis}
  Suppose that $X$ and $Y$ are analytic metric spaces and $\phi
  \from X \to Y$ is a \Borel function that is not \Baire class one. Then
  there exist $\phi_{\Bairespace} \in \set{\constant{\Bairespace}}
  \union \set{\inclusion{\Bairespace}{Z}}[Z \in \set{\Bairespace,
  \compactifiedBairespace}]$ and $\phi_{\compactifiedBairespace
  \setminus \Bairespace} \in \set{\constant{\compactifiedBairespace
  \setminus \Bairespace}} \union \set{\inclusion
  {\compactifiedBairespace \setminus \Bairespace}{Z}}[Z \in \set
  {\compactifiedBairespace \setminus \Bairespace,
  \compactifiedBairespace}] \union \set{\inclusion{\Bairetree}{Z}
  \composition \discrete}[Z \in \set{\Bairetree,
  \onepointcompactification{\Bairetree}, \compactifiedextendedBairespace \setminus \Bairespace,
  \extendedBairespace, \compactifiedextendedBairespace}]$ for
  which there is a closed continuous embedding of
  $\phi_{\Bairespace} \disjointunion \phi_{\compactifiedBairespace
  \setminus \Bairespace}$ into $\phi$.
\end{theorem}

\begin{theoremproof}
  \Hurewicz's dichotomy theorem for \Fsigma sets yields a closed
  continuous embedding $\psi \from \compactifiedBairespace \to X$
  with $\image{(\phi \composition \psi)}{\Bairespace} \intersection
  \closure{\image{(\phi \composition \psi)}{\compactifiedBairespace
  \setminus \Bairespace}} = \emptyset$. As $\pair{\psi}{\identity{\closure
  {\image{(\phi \composition \psi)}{\compactifiedBairespace}}}}$ is a
  closed continuous embedding of $\phi \composition \psi$ into
  $\phi$, by replacing the latter with the former, we can assume that
  $X = \compactifiedBairespace$ and $\image{\phi}{\Bairespace}
  \intersection \closure{\image{\phi}{\compactifiedBairespace
  \setminus \Bairespace}} = \emptyset$.
  
  By Proposition \ref{Bairespace:continuous}, there is a
  $\meet$-embedding $\pi \from \Bairetree \to \Bairetree$ for which
  $\restriction{(\phi \composition \extension{\pi})}{\Bairespace}$ is
  continuous. By composing $\pi$ with the $\meet$-embedding
  given by Proposition \ref{Baireclassone:disjoint}, we can assume
  that $\closure{\image{(\phi \composition \extension{\pi})}
  {\Bairespace}} \intersection \closure{\image{(\phi \composition
  \extension{\pi})}{\compactifiedBairespace \setminus \Bairespace}}
  = \emptyset$. By composing $\pi$ with the $\meet$-embedding
  given by Theorem \ref{Bairespace:closedcontinuousembedding},
  we can assume that $\restriction{(\phi \composition \extension{\pi})}
  {\Bairespace}$ is constant or extends to a closed continuous
  embedding on $\Bairespace$ or $\compactifiedBairespace$. And
  by composing $\pi$ with the $\meet$-embedding given by Theorem
  \ref{Bairetree:closedcontinuousembedding}, we can assume that
  $\restriction{(\phi \composition \extension{\pi})}
  {\compactifiedBairespace \setminus \Bairespace}$ is constant,
  $\restriction{(\phi \composition \extension{\pi})}
  {\compactifiedBairespace \setminus \Bairespace}$ extends to a
  closed continuous embedding on $\compactifiedBairespace
  \setminus \Bairespace$ or $\compactifiedBairespace$, or $\phi
  \composition \extension{\pi} \composition \inverse{\discrete}$
  extends to a closed continuous embedding on $\Bairetree$,
  $\onepointcompactification{\Bairetree}$, $\compactifiedextendedBairespace \setminus \Bairespace$,
  $\extendedBairespace$, or $\compactifiedextendedBairespace$.

  By Proposition \ref{Bairespace:basis}, there exist
  $\phi_{\Bairespace} \in \set{\constant{\Bairespace}} \union \set
  {\inclusion{\Bairespace}{Z}}[Z \in \set{\Bairespace,
  \compactifiedBairespace}]$ and $\psi_{\Bairespace} \from \closure
  {\image{\phi_{\Bairespace}}{\Bairespace}} \to \closure{\image{\phi}
  {\Bairespace}}$ for which $\pair{\restriction{\extension{\pi}}
  {\Bairespace}[\Bairespace]}{\psi_{\Bairespace}}$ is a closed
  continuous embedding of $\phi_{\Bairespace}$ into $\restriction
  {\phi}{\Bairespace}$. By Proposition \ref{Bairetree:basis}, there
  exist $\phi_{\compactifiedBairespace \setminus \Bairespace} \in
  \set{\constant{\compactifiedBairespace \setminus \Bairespace}}
  \union \set{\inclusion{\compactifiedBairespace \setminus
  \Bairespace}{Z}}[Z \in \set{\compactifiedBairespace \setminus
  \Bairespace, \compactifiedBairespace}] \union \set{\inclusion
  {\Bairetree}{Z} \composition \discrete}[Z \in \set{\Bairetree,
  \onepointcompactification{\Bairetree}, \compactifiedextendedBairespace \setminus \Bairespace,
  \extendedBairespace, \compactifiedextendedBairespace}]$ and
  $\psi_{\compactifiedBairespace \setminus \Bairespace} \from
  \closure{\image{\phi_{\compactifiedBairespace \setminus
  \Bairespace}}{\compactifiedBairespace \setminus \Bairespace}}
  \to \closure{\image{\phi}{\compactifiedBairespace \setminus
  \Bairespace}}$ for which $\pair{\restriction{\extension{\pi}}
  {\compactifiedBairespace \setminus \Bairespace}
  [\compactifiedBairespace \setminus \Bairespace]}
  {\psi_{\compactifiedBairespace \setminus \Bairespace}}$ is a
  closed continuous embedding of $\phi_{\compactifiedBairespace
  \setminus \Bairespace}$ into $\restriction{\phi}
  {\compactifiedBairespace \setminus \Bairespace}$. Then $\pair
  {\restriction{\extension{\pi}}{\compactifiedBairespace}
  [\compactifiedBairespace]}{\psi_{\Bairespace} \disjointunion
  \psi_{\compactifiedBairespace \setminus \Bairespace}}$ is a
  closed continuous embedding of $\phi_{\Bairespace} \disjointunion
  \phi_{\compactifiedBairespace \setminus \Bairespace}$ into $\phi$.
\end{theoremproof}

Theorems \ref{sigmacontinuous:basis} and \ref{Baireclassone:basis}
together provide the promised twenty-seven element basis under
closed continuous embeddability for the class of
non-$\sigma$-continuous-with-closed-witnesses \Borel functions
between analytic metric spaces.

\bibliographystyle{amsalpha}
\bibliography{bibliography}

\end{document}